# Exact semi-separation of variables in waveguides with nonplanar boundaries

G.A. Athanassoulis ([1]) ([2]), Ch.E. Papoutsellis ([3])

**Abstract**

Series expansions of unknown fields $\Phi = \sum \varphi_n Z_n$ in elongated waveguides are commonly used in acoustics, optics, geophysics, water waves and other applications, in the context of coupled-mode theories (CMTs). The transverse functions $Z_n$ are determined by solving local Sturm-Liouville problems (reference waveguides). In most cases, the boundary conditions assigned to $Z_n$ cannot be compatible with the physical boundary conditions of $\Phi$, leading to slowly convergent series, and rendering CMTs mild-slope approximations. In the present paper, the heuristic approach introduced in (Athanassoulis & Belibassakis 1999 *J. Fluid Mech.* **389**, 275-301) is generalized and justified. It is proved that an appropriately enhanced series expansion becomes an exact, rapidly-convergent representation of the field $\Phi$, valid for any smooth, nonplanar boundaries and any smooth enough $\Phi$. This series expansion can be differentiated termwise everywhere in the domain, including the boundaries, implementing an exact semi-separation of variables for non-separable domains. The efficiency of the method is illustrated by solving a boundary value problem for the Laplace equation, and computing the corresponding Dirichlet-to-Neumann operator, involved in Hamiltonian equations for nonlinear water waves. The present method provides accurate results with only a few modes for quite general domains. Extensions to general waveguides are also discussed.

**Table of contents**

---

([1]) Professor, National Technical University of Athens, School of Naval Architecture and Marine Engineering, Athens, Greece, Email: mathan@central.ntua.gr, makathan@gmail.com

([2]) Research Professor, Research Center for High Performance Computing, ITMO University, St. Petersburg, Russia

([3]) PhD student, National Technical University of Athens, School of Naval Architecture and Marine Engineering, Athens, Greece, Email: cpapoutse@central.ntua.gr, cpapoutsellis@gmail.com

# 1. Introduction

In many applications we need to solve boundary-value problems (BVPs) in non-uniform, strip-like domains. Prominent examples are the propagation of surface gravity waves or sound waves in nearshore sea environments, acoustic waves in nonuniform ducts, electromagnetic waves in curved or varying cross-section microwave and optical waveguides, etc. In many such cases, the domain filled by the medium, denoted by $D_h^\eta$, is contained between a lower boundary $\Gamma_h : z = -h(\boldsymbol{x})$ and an upper boundary $\Gamma^\eta : z = \eta(\boldsymbol{x})$ ([4]), where $\boldsymbol{x} = (x_1, x_2)$ denote the horizontal coordinates and $z$ is the vertical (transverse) coordinate. The domain $D_h^\eta$ may be laterally unbounded, semi-bounded, as in Fig. 1, or completely bounded. Very often, in these problems, the horizontal extent is much larger than the depth, and the domain contains a large number of wavelengths, which makes the implementation of direct numerical simulation (e.g., 3D finite-element method) very demanding and sometimes practically impossible. On the other hand, if the domain is a uniform strip, then usually the method of separation of variables provides us with an analytical or semi-analytical solution in the form of eigenfunctions expansions. This solution technique has been extended to domains with slowly varying boundaries in the context of various *coupled-mode theories* (CMTs), by approximating the unknown field $\Phi(\boldsymbol{x}, z)$ as a *local-mode series expansion* of the form

$$\Phi(\boldsymbol{x}, z) = \sum_n \varphi_n(\boldsymbol{x}) Z_n(z; \boldsymbol{x}), \tag{1.1}$$

where $Z_n$ are known functions and $\varphi_n$ are unknown modal amplitudes. As a rule, the *vertical functions* $Z_n$ are constructed by means of a *localized version* of the transverse (vertical) Sturm-Liouville problem, obtained by separation of variables in a uniform strip adapted to the interval $-h(\boldsymbol{x}) \leq z \leq \eta(\boldsymbol{x})$ (the *reference waveguide at* $\boldsymbol{x}$, [1], Sec. 7.1). Using Eq. (1.1), it is possible to approximately reduce the wave propagation problem in $D_h^\eta$ to a *coupled-mode system* of horizontal differential equations, implementing a dimension reduction of the initial problem. The CMT offers a promising solution method in cases where a small number of modes provides a satisfactory approximation to the unknown field $\Phi(\boldsymbol{x}, z)$.

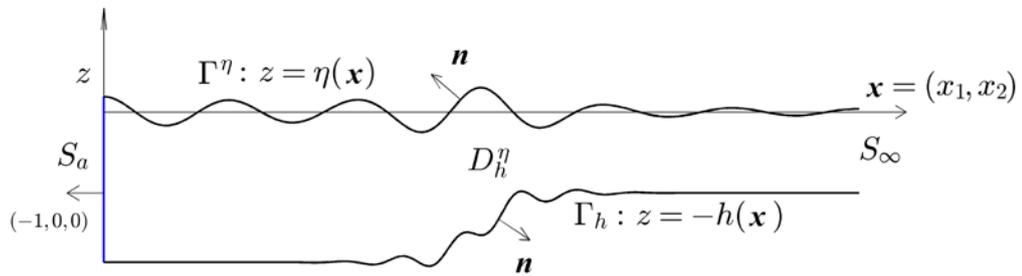

**Figure 1**. Geometrical configuration of the generic, semi-infinite waveguide

---

([4]) We shall often be using the terminology *seabed* (or *bottom*) for the lower boundary, and *free surface* for the upper boundary, corresponding to the case where $D_h^\eta$ is an oceanic environment, although the method and the results presented herein can be used in other applications as well.

The origin of the CMT goes back to early 50s, first appearing in connection with electromagnetic [2], and acoustic [3] waveguides of varying cross sections (horns). Since then, the CMT has been applied to the study of almost all types of elongated waveguides([5]): microwave and optical waveguides (sometimes under the name *generalized telegraphist's equations*) [4], [5], [6]; duct acoustics (horns) [7], [8], [9]; duct acoustics with mean flow [10]; hydroacoustics [1] (Sec. 7.1), [11], [12], [13], [14]; seismology and geophysics [15], [16]; water waves [17], [18], [19], [20], [21], [22], [23], [24], [25], [26], [27], [28]; hydroelasticity [23].

The local vertical eigenfunctions $Z_n$, as obtained by the reference waveguide, satisfy boundary conditions which, in general, are incompatible with the boundary conditions of the field $\Phi(x,z)$. As a result, the series representation (1.1) is not able to reproduce correctly the boundary behavior of $\Phi(x,z)$ on $z = -h(x)$ and $z = \eta(x)$, except for the cases of a Dirichlet condition and a homogeneous Neumann condition at horizontal boundaries. For this reason, the classical CMTs are considered as mild-slope approximations. This deficiency is well known and various remedies have been proposed in the literature. Many authors have tried to avoid the problem by replacing the smoothly varying, nonplanar boundary by a piecewise planar (staircase) boundary; see e.g. [29], [30] for applications to acoustic problems, [31], [32] for water-wave problems, and [33] for electromagnetic waveguides. This method may provide satisfactory solutions, at points far from the changed boundary, for 2D and axisymmetric geometries, but it is inappropriate for fully 3D environments and for complicated boundary conditions. Rutherford and Hawker [34] used perturbation theory to construct a first-order (with respect to the slop $\nabla_x h$) correction $\delta Z_n$ to each vertical function $Z_n$, in the context of hydroacoustics. In late 90's, Athanassoulis and Belibassakis, in the context of water-wave theory [19], [20], [21], proposed another solution of this problem, by introducing additional modes in the series expansion (1.1). This approach, initially introduced in a heuristic, ad hoc way, turned to be very efficient. The enhanced series expansion exhibits a faster convergence rate ($O(n^{-4})$ instead of $O(n^{-2})$), and its use permitted to obtain satisfactory numerical solutions to problems concerning waveguides with strongly varying boundaries. The CMT arising by using the enhanced coupled-mode series, called *consistent coupled-mode theory* (CCMT) by its inventors, developed further by the same and other authors; see e.g. [35], [23], [36], [8], [9], [37]. Although the enhanced modal expansion has been extensively used for more than fifteen years, the mathematical reasons of its success are not well understood. This paper is mainly devoted to the study of mathematical aspects related to the following two questions:

Q1. Is it possible to construct a convergent series expansion, like (1.1), for any field $\Phi(x,z)$ defined in a strip-like domain $D_h^\eta$ with arbitrary boundaries $\Gamma_h : z = -h(x)$ and $\Gamma^\eta : z = \eta(x)$ ? And, if the answer is *yes*,

Q2. Under what conditions this series expansion can be differentiated term-by-term, providing convergent series expansions to the first and second derivatives of $\Phi(x,z)$ ?

The answer to question Q1 is yes, provided that the field $\Phi$ and the boundaries $\Gamma_h$, $\Gamma^\eta$ are smooth, and the series expansion includes all eigenfunctions of the reference waveguide (Sturm-Liouville problem), say $\{Z_n, n \geq 0\}$, plus two additional modes, $Z_{-2}$, $Z_{-1}$, appro-priately defined so that to remove any restrictions of the series on the boundaries $\Gamma^\eta$ and $\Gamma_h$ (Theorem 1, in section 2). Special cases of this result have been proved in [35], [38]. Question Q2 is considered here for the first time. The answer is that the termwise differentiated series are convergent to the corresponding derivatives of the field $\Phi$ under additional smoothness assumptions (Theorem 2, in section 4).

---

([5]) The list of application areas of CMT is not exhaustive, and the cited literature is just indicative.

These results establish theoretically the validity of the semi-separation of variables in non-separable geometries (strip-like domains with nonplanar boundaries), and ensure that the CCMT provides an *exact reformulation* of BVPs in these geometries. The numerical performance of the proposed method is also addressed in this paper (section 5), by studying simple benchmark problems with known analytical solutions. An interesting by-product of the present solution method is the highly accurate computation of a Dirichlet-to-Neumann (DtN) operator which is of fundamental importance for solving fully nonlinear water-wave problems.

Recalling that the set of eigenfunctions $\{Z_n, n \geq 0\}$ of the local vertical Sturm-Liouville problem (reference waveguide) is an $L^2-$basis in the vertical interval $(-h, \eta)$, the addition of two more elements, $Z_{-2}$, $Z_{-1}$, raises the question of linear independence. To answer this question correctly, it is essential to pose it in the appropriate function space: the system $\{Z_{-2}, Z_{-1}, Z_n, n \geq 0\}$ is not linearly independent in the space $L^2(-h, \eta)$, but *it is linearly independent* in the Sobolev space $H^2(-h, \eta)$, which is the appropriate space for problems in which the trace of the first derivative on the boundaries should be well defined. Such kind of questions have also been discussed, in conjunction with the application of variational methods to some simple problems, by Storch and Strang [39]; see also [40].

The structure of the remaining part of this paper is as follows. In section 2 we show that the enhanced local-mode series converges rapidly in $D_h^\eta$, independently of the shape of the (smooth) boundaries $\Gamma^\eta$ and $\Gamma_h$, and the boundary behaviour of the (smooth) field $\Phi$. In section 3 we derive analytical expressions for the first and second derivatives $\partial_{x_i} Z_n$, $\partial_{x_i}^2 Z_n$, $i = 1, 2$, of the vertical functions, and obtain estimates of these quantities with respect to $n$ (as $n \to \infty$). These results are used, in section 4, to derive a priori decay rates for the derivatives of the modal amplitudes, $\partial_{x_i} \varphi_n$ and $\partial_{x_i}^2 \varphi_n$, with respect to $n$, which are then exploited to prove that the enhanced local-mode series can be differentiated term-by-term twice, under appropriate smoothness assumptions. In section 5(a) we apply the proposed semi-separation of variables to the numerical solution of specific BVPs for the Laplace equation, with non-flat upper boundary. By selecting cases for which analytical solutions are known, we show that the first few modes converge to their exact values as $N_{tot}^{-6.5}$, and this rate gradually deteriorates to $N_{tot}^{-3.5}$ for the last retained mode, where $N_{tot}$ is the total number of modes considered in the truncated expansion. As a result, the field $\Phi$ converges to its exact value as $N_{tot}^{-3.5}$, while the corresponding DtN operator, which is expressed in terms of $\varphi_{-2}$ and the data of the problem, inherits the "superconvergence" rate $N_{tot}^{-6.5}$. This fact is expected to be of great value for solving fully nonlinear water-wave problems, since only the DtN operator (and not the whole field $\Phi$) is involved in the corresponding nonlinear evolution equations [41], [27]. Extensions of the results obtained in this paper to more general reference waveguides are dis-cussed in section 6. Various technical issues and extended calculations are given in detail in five appendices.

## 2. Series expansion of a smooth function defined throughout a strip with non-planar boundaries

Consider, for definiteness, the nonplanar, semi-infinite, strip-like domain

$$D_h^\eta(X) = \{(x,z) \in X \times \mathbb{R} : x \in X, -h(x) < z < \eta(x)\},$$

where $X = \{x = (x_1, x_2) \in \mathbb{R}^2 : x_1 \geq a\}$ is the common projection of the boundary surfaces $\Gamma_h$ and $\Gamma^\eta$ on the $x = (x_1, x_2)$ plane. The lateral vertical boundary at $x_1 = a$ is denoted by $S_a$, and the remaining, infinite, part of the lateral boundary is denoted by $S_\infty$; see Fig. 1. The local depth (height) of the non-uniform strip is denoted by $H(x)$, and we assume that

$$H = H(x) = \eta(x) + h(x) > 0. \tag{2.1}$$

For later use we also define the subdomains $D_h^\eta(\tilde{X}_{part})$ consisting of all points of $D_h^\eta(X)$ with $x \in \tilde{X}_{part}$, $\tilde{X}_{part}$ being an open subregion, finite or infinite, of $X$. A scalar field $\Phi = \Phi(x,z)$ is defined on the closed domain $\overline{D}_h^\eta(X)$, where $\overline{D}_h^\eta(X) = D_h^\eta(X) \cup \Gamma_h \cup \Gamma^\eta \cup S_a \cup S_\infty$. Smoothness conditions on the field $\Phi$ and the boundaries $\Gamma_h$, $\Gamma^\eta$ will be introduced subsequently. In this section we shall prove that $\Phi$ can be expanded in a rapidly convergent series of the form

$$\Phi(x,z) = \varphi_{-2}(x) Z_{-2}(z;\eta(x),h(x)) + \varphi_{-1}(x) Z_{-1}(z;\eta(x),h(x))$$
$$+ \sum_{n=0}^{\infty} \varphi_n(x) Z_n(z;\eta(x),h(x)), \tag{2.2}$$

where $Z_n$, $n \geq 0$ are the eigenfunctions of a reference waveguide, and the additional functions $Z_{-2}$, $Z_{-1}$ are appropriately selected boundary modes.

(a) *The reference waveguide and the boundary modes*

In the present sybsection, the *local vertical Sturm-Liouville problem* (equivalently, the *reference waveguide*) is chosen as follows:

$$\partial_z^2 Z_n(z) + \left(k_n^2 - Q\right) Z_n(z) = 0, \quad -h(x) < z < \eta(x), \tag{2.3a}$$

$$B^\eta Z_n \equiv [(\partial_z - \mu_0) Z_n]_{z=\eta} = 0, \, (^6) \tag{2.3b}$$

$$B_h Z_n \equiv [\partial_z Z_n]_{z=-h} = 0, \tag{2.3c}$$

---

(⁶) The notation $[G(x,z)]_{z=\eta}$ is used to denote the trace of function $G$ on the boundary $\Gamma^\eta$. Similar notation is used for the traces on the other boundaries.

where $Q$, $\mu_0$ are real constants. All the results obtained in this section remain valid for the fairly general Sturm-Liouville problem with $Q = Q(x,z)$ and general, first-order boundary operators; see Remark 2.2.

The eigensystem $\{k_n, Z_n(z), n \geq 0\}$ is dependent on $\mu_0$, $Q$, $\eta = \eta(x)$ and $h = h(x)$. Thus, $k_n = k_n(x)$ and $Z_n(z) = Z_n(z;x)$ ([7]). In the present work, the eigenfunctions $Z_n(z)$ are normalized so that

$$Z_n(z = \eta(x)) = 1. \tag{2.4}$$

From the general theory of Sturm-Liouville problems we recall the following facts, which are needed in the proof of the convergence of the series expansion (2.2): **i)** the system of eigenfunctions $\{Z_n(z), n \geq 0\}$ forms an orthogonal basis of the space of functions $L^2(-h, \eta)$, and **ii)** for large values of $n$, the quantities $k_n$, $Z_n(z)$ and $\|Z_n\|_{L^2}$ exhibit the following asymptotic behavior:

$$k_n = \frac{\pi n}{H(x)} + \frac{O(1)}{n} = O(n), \qquad Z_n = \frac{\cos(k_n(z + h(x)))}{\cos(k_n H(x))} + \frac{O(1)}{n} \tag{2.5a,b}$$

$$\|Z_n\|_{L^2}^2 \equiv \int_{-h(x)}^{\eta(x)} |Z_n(z)|^2 dz = \frac{H(x)}{2} + \frac{O(1)}{n}. \tag{2.5c}$$

Asymptotic results (2.5) are classical and can be found in various books; see, e.g., [42], Chapter 1, [43], Chapter 10. It is worth noting that the leading-order asymptotic terms, shown in Eqs. (2.5), are not dependent on $Q$; they are only dependent on the boundary functions $\eta = \eta(x)$ and $h = h(x)$, revealing the critical role played by the latter on the asymptotics and, through them, on the rate of convergence of the expansion (2.2), as we shall see in the sequel.

The boundary modes $Z_{-2}$ (associated with $\Gamma^\eta$) and $Z_{-1}$ (associated with $\Gamma_h$) are defined so that the enhanced series expansion (2.2) is able to recover directly the values of $B^\eta \Phi \neq 0$ and $B_h \Phi \neq 0$ (compare with Eqs. (2.3b) and (2.3c)). Thus, the only essential conditions they have to satisfy are $B^\eta Z_{-2} \neq 0$ and $B_h Z_{-1} \neq 0$. In this connection, $Z_{-2}$ and $Z_{-1}$, are selected to be smooth functions satisfying the following boundary conditions

$$[Z_{-2}]_{z=\eta} = 1, \quad B^\eta Z_{-2} = 1/h_0, \quad B_h Z_{-2} = 0, \tag{2.6a,b,c}$$

$$[Z_{-1}]_{z=\eta} = 1, \quad B^\eta Z_{-1} = 0, \quad B_h Z_{-1} = 1/h_0, \tag{2.7a,b,c}$$

---

([7]) The local vertical functions $Z_n = Z_n(z)$, $n \geq -2$, will be denoted by anyone of the symbols $Z_n$, $Z_n(z)$, $Z_n(z;\eta,h)$, $Z_n(z;\eta(x),h(x))$ or $Z_n(x,z)$, as appropriate in the specific context, in order to make clear an argument, or a treatment, or an implication. For the same reason, various alternative expressions will also be used for denoting the local eigenvalues $k_n = k_n(H) = k_n(\eta,h) = k_n(\eta(x),h(x)) = k_n(x)$.

where the constant $h_0$ (reference depth) is introduced for dimensional consistency. Conditions (2.6a), (2.7a) serve for normalization (compare with Eq. (2.4)), while conditions (2.6c), (2.7b) are imposed in order to separate the influence of the two boundary modes $Z_{-2}$, $Z_{-1}$; the former will act on the upper surface $\Gamma^\eta$ but not on the lower surface $\Gamma_h$, and the latter inversely. Clearly, Eqs. (2.6), (2.7) do not define uniquely the functions $Z_{-2}$, $Z_{-1}$. Among the infinite possible choices, we pick out the one corresponding to the least degree polynomials (in the vertical variable $z$) ([8]). Then, we easily find

$$Z_{-2}(z;\eta,h) = \frac{\mu_0 h_0 + 1}{2 h_0} \frac{(z+h)^2}{(\eta+h)} - \frac{\mu_0 h_0 + 1}{2 h_0}(\eta+h) + 1, \qquad (2.8a)$$

$$Z_{-1}(z;\eta,h) = \frac{\mu_0 h_0 - 1}{2 h_0} \frac{(z+h)^2}{(\eta+h)} + \frac{1}{h_0}(z+h) - \frac{\mu_0 h_0 + 1}{2 h_0}(\eta+h) + 1. \qquad (2.8b)$$

(b) *The fundamental expansion theorem*

The section $z \to F(\mathbf{x}, z)$ of a function $F: \bar{D}_h^\eta(X) \to \mathbb{R}$, will be denoted as $F(\mathbf{x}, \cdot)$.

**Definition 1**. We shall say that $F(\mathbf{x}, \cdot)$ *satisfies a local vertical smoothness condition* $VSC(k; \mathbf{x})$ if $F(\mathbf{x}, \cdot) \in H^k(-h, \eta)$. If, further, the norm $\|F(\mathbf{x}, \cdot)\|_{H^k}$ is uniformly bounded in $\tilde{X}_{part}$, then, we shall say that $F(\mathbf{x}, \cdot)$ *satisfies a uniform vertical smoothness condition* $VSC(k; \tilde{X}_{part})$. $\tilde{X}_{part}$ may be replaced by $X$, if the uniform boundedness of the Sobolev norm holds true in $X$. □

$H^k(-h, \eta)$ is the standard Sobolev space of functions defined on $[-h, \eta]$; see e.g. [44] Chapter 8. Recall that $H^k(-h, \eta)$ is continuously embedded in $C^{k-1}([-h, \eta])$ ([44], p. 217), which ensures that the boundary values of the vertical derivatives $\partial_z^m F(\mathbf{x}, z, t)$, at $z = \eta, -h$, are well-defined and bounded for $0 \leq m \leq k-1$.

**Theorem 1 [The fundamental expansion theorem].** *Assume that $\Phi(\mathbf{x}, \cdot)$ satisfies a local $VSC(k=2; \mathbf{x})$. Then, the quantities*

$$\varphi_{-2}(\mathbf{x}) = h_0 B^\eta \Phi, \quad \varphi_{-1}(\mathbf{x}) = h_0 B_h \Phi, \qquad (2.9a,b)$$

*and*

$$\Phi^*(\mathbf{x}, z) = \Phi(\mathbf{x}, z) - \varphi_{-2}(\mathbf{x}) Z_{-2} - \varphi_{-1}(\mathbf{x}) Z_{-1}, \qquad (2.9c)$$

*are well-defined, and*

(i) *The series expansion (2.2) converges uniformly on $[-h(\mathbf{x}), \eta(\mathbf{x})]$, and the modal amplitudes $\varphi_n(\mathbf{x})$, $n \geq 0$, are expressed by the equations*

---

([8]) Other choices are possible without affecting the validity of the expansion.

$$\varphi_n(\boldsymbol{x}) = \|Z_n\|_{L^2}^{-2} \int_{-h(\boldsymbol{x})}^{\eta(\boldsymbol{x})} \Phi^* Z_n \, dz, \quad n \geq 0. \tag{2.9d}$$

**(ii)** *The modal amplitudes satisfy the following local asymptotic estimate* (*with respect to* $n$)

$$|\varphi_n(\boldsymbol{x})| \leq A_2(\boldsymbol{x}) \, n^{-2}. \tag{2.10}$$

**(iii)** *If, further,* $\Phi(\boldsymbol{x}, \cdot)$ *satisfies the local VSC* ($k = 4; \boldsymbol{x}$), *then the asymptotic estimate* (2.10) *is improved to*

$$|\varphi_n(\boldsymbol{x})| \leq A_4(\boldsymbol{x}) \, n^{-4}. \tag{2.11}$$

*The latter estimate cannot be improved further under the present choice of the boundary modes* $Z_{-2}, Z_{-1}$.

**Proof**: The modal amplitudes $\varphi_{-2}$ and $\varphi_{-1}$ are well-defined since $\Phi(\boldsymbol{x}, \cdot) \in H^2(-h, \eta)$. Accordingly, $\Phi^*(\boldsymbol{x}, z)$ is also well-defined.

**(i)** Applying the boundary operators $B^\eta$ and $B_h$ to both members of Eq. (2.9c) and taking into account the boundary conditions of $Z_{-2}$ and $Z_{-1}$, Eqs. (2.6), (2.7), we readily find that $B_h \Phi^* = 0$ and $B^\eta \Phi^* = 0$. Thus, $\Phi^*(\boldsymbol{x}, z)$, as a function of the vertical coordinate $z$, satisfies the same boundary conditions as the eigenfunctions $Z_n(z)$, $n \geq 0$, of the regular Sturm-Liouville problem (2.3). As a consequence, the $L^2$–series expansion $\Phi^*(\boldsymbol{x}, z) = \sum_{n=0}^{n=\infty} \varphi_n(\boldsymbol{x}) Z_n(z)$ is also *uniformly convergent* (see Theorem 4.1, p. 197, of [45]). Eq. (2.9d) is a direct consequence of the orthogonality of the vertical eigenfunctions. Substituting the series expansion of $\Phi^*$ in Eq. (2.9c), and solving for $\Phi(\boldsymbol{x}, z)$, we obtain the series expansion (2.2), having already established its uniform convergence.

**(ii)** We proceed now to prove the estimate (2.10). Using Eq. (2.3a), we reformulate Eq. (2.9d) as

$$\varphi_n(\boldsymbol{x}) = \|Z_n\|_{L^2}^{-2} k_n^{-2} Q_n^* \int_{-h}^{\eta} \Phi^* \partial_z^2 Z_n \, dz, \tag{2.12}$$

where

$$Q_n^* \equiv -\frac{1}{1 - Q/k_n^2} = O(1), \tag{2.13}$$

in view of the estimate (2.5a). Under the assumption $\Phi^*(\boldsymbol{x}, \cdot) \in H^2(-h, \eta)$, it is legitimate to perform two integrations by parts in Eq. (2.12), obtaining

$$\varphi_n(\boldsymbol{x}) = \|Z_n\|_{L^2}^{-2} k_n^{-2} Q_n^* \left( \left[ \Phi^* \partial_z Z_n \right]_{-h}^{\eta} - \left[ \partial_z \Phi^* \, Z_n \right]_{-h}^{\eta} + \int_{-h}^{\eta} \partial_z^2 \Phi^* \, Z_n \, dz \right). \tag{2.14}$$

The boundary terms in the parenthesis, in the right-hand side of Eq. (2.14), cancel out since $\partial_z \Phi^* = \partial_z Z_n = 0$ at $z = -h$, and $\Phi^* \partial_z Z_n - \partial_z \Phi^* Z_n = \Phi^* \mu_0 - \partial_z \Phi^* = -B^\eta \Phi^* = 0$, at $z = \eta$. Hence,

$$\varphi_n(\boldsymbol{x}) = \|Z_n\|_{L^2}^{-2} k_n^{-2} Q_n^* \int_{-h}^{\eta} \partial_z^2 \Phi^* \, Z_n \, dz, \tag{2.15}$$

and the estimate (2.10) follows after a straightforward application of Cauchy-Swartz inequality in conjunction with the estimates (2.5) and (2.13).

**(iii)** In Eq. (2.15) we substitute $Z_n$, in the integrand, by $k_n^{-2} Q_n^* \partial_z^2 Z_n$ (from Eq. (2.3a)), obtaining

$$\varphi_n(\boldsymbol{x}) = \|Z_n\|_{L^2}^{-2} k_n^{-4} (Q_n^*)^2 \left( \int_{-h}^{\eta} \partial_z^2 \Phi^* \, \partial_z^2 Z_n \, dz \right). \tag{2.16}$$

Under the assumption $\Phi^*(\boldsymbol{x},\cdot) \in H^4(-h,\eta)$, it is legitimate to perform two more integrations by parts, obtaining

$$\varphi_n(\boldsymbol{x}) = \|Z_n\|_{L^2}^{-2} k_n^{-4} (Q_n^*)^2 \left( \left[ \partial_z^2 \Phi^* \partial_z Z_n \right]_{-h}^{\eta} - \left[ \partial_z^3 \Phi^* Z_n \right]_{-h}^{\eta} + \int_{-h}^{\eta} \partial_z^4 \Phi^* Z_n \, dz \right). \tag{2.17}$$

The boundary terms in the parentheses are all bounded ($Z_n(\eta) = 1$ due to Eq. (2.4), $\partial_z Z_n(\eta) = \mu_0 Z_n(\eta) = \mu_0$ due to Eqs. (2.3b) and (2.4), $\partial_z Z_n(-h) = 0$ due to Eq. (2.3c); for estimating $Z_n(-h)$ use is made of Eq. (2.5b), obtaining $|Z_n(-h)| \leq |\cos(k_n H)|^{-1} + O(n^{-1}) \leq 2$ for $n$ greater than some fixed $n_*$; the derivatives $\partial_z^2 \Phi^*$ and $\partial_z^3 \Phi^*$ are well-defined since $\Phi^*(\boldsymbol{x},\cdot) \in H^4(-h,\eta)$ ). Further, the integral is also bounded as can be seen by applying the Cauchy-Swartz inequality. Hence, Eq. (2.17), in conjunction with the estimates (2.5a) and (2.13), leads to the estimate (2.11).

Note that, this estimate cannot be improved further because of the presence of the boundary terms in the right-hand side of Eq. (2.17). These boundary terms cannot be zero, as far as $\partial_z^2 \Phi^*$, $\partial_z^3 \Phi^*$ are free of any restrictions. Thus, even if $\Phi^*$ is smoother, e.g. $\Phi^*(\boldsymbol{x},\cdot) \in H^6(-h,\eta)$, in which case the last (integral) term in the right-hand side of Eq. (2.17) can be reformulated, by using again the substitution $Z_n = k_n^{-2} \partial_z^2 Z_n$, and becomes $O(k_n^{-2})$, the presence of the boundary terms containing $\partial_z^2 \Phi^*$, $\partial_z^2 \Phi^*$ blocks the order of $\varphi_n(\boldsymbol{x})$ to $O(k_n^{-4})$. ∎

**Corollary 1.** *If $\Phi(\boldsymbol{x},\cdot)$ satisfies the uniform VSC($k=2;\tilde{X}_{\text{part}}$) (resp. VSC($k=4;\tilde{X}_{\text{part}}$)), then the series expansion* (2.2) *is uniformly convergent in the subdomain $D_h^\eta(\tilde{X}_{\text{part}})$, and the estimate* (2.10) *(resp.* (2.11)*) is uniformly valid in $\tilde{X}_{\text{part}}$, that is $\|\varphi_n\|_{\infty,\tilde{X}_{\text{part}}} = O(n^{-2})$ (resp. $\|\varphi_n\|_{\infty,\tilde{X}_{\text{part}}} = O(n^{-4})$). In all the above statements $\tilde{X}_{\text{part}}$ can be replaced by $X$, if the assumptions hold true globally in $X$.* □

**Remarks: 2.1.** Theorem 1 and a sketch of its proof were announced in 2003 [35] for a flat upper boundary. A proof of assertion (i) for the case $Q = 0$, and more restrictive assumptions on the field $\Phi$ and boundaries $\Gamma^\eta$ and $\Gamma_h$, was presented in 2006 [38].

**2.2.** Asymptotic estimates (2.5) and Theorem 1 remain valid for the fairly general Sturm-Liouville problem

$$\partial_z^2 Z_n(z) + \left( k_n^2 - Q(\boldsymbol{x},z) \right) Z_n(z) = 0, \quad -h(\boldsymbol{x}) < z < \eta(\boldsymbol{x}), \tag{2.18a}$$

$$B^\eta Z_n \equiv [(\partial_z - \mu_0(\pmb{x}))Z_n]_{z=\eta} = 0, \qquad (2.18b)$$

$$B_h Z_n \equiv [(\partial_z - \nu_0(\pmb{x}))Z_n]_{z=-h} = 0, \qquad (2.18c)$$

under the assumption that $\mu_0(\pmb{x})$, $\nu_0(\pmb{x})$ and $Q(\pmb{x},z)$ are continuous functions of $\pmb{x}$, and $Q(\pmb{x},z)$ is sufficiently smooth with respect to $z$. The proof of Theorem 1 in that case requires additionally the estimation of the $z$–derivatives of $Q_n^* = Q_n^*(\pmb{x},z)$ and, thus, becomes more technical.

**2.3.** The validity of Theorem 1 and its corollary is not affected by the choice of boundary modes $Z_{-2}$, $Z_{-1}$ as far as the conditions (2.6) and (2.7) hold true.

**2.4.** By the procedure described in this section, we have constructed a Riesz basis of the Sobolev space $H^2(-h,\eta)$. In fact, it can be proved that, if the system $\{Z_n, n \geq 0\}$ is an orthonormal basis of $L^2(-h,\eta)$, then, the system $\{Z_{-2}, Z_{-1}, Z_n/k_n^2, n \geq 0\}$ is a Riesz basis of $H^2(-h,\eta)$. A direct proof of this result will be published elsewhere.

## 3. The $x$ – derivatives of the vertical eigensystem

In order to study the $x$–differentiability properties of the infinite series expansion (2.2), we need to calculate the $x$–derivatives of the vertical eigenfunctions $Z_n(z;\pmb{x}) = Z_n(z;\eta(\pmb{x}),h(\pmb{x}))$, and determine their asymptotic behavior for large values of $n$. Since

$$\partial_{x_i} Z_n(z;\pmb{x}) = \partial_\eta Z_n(z;\eta,h)\,\partial_{x_i}\eta(\pmb{x}) + \partial_h Z_n(z;\eta,h)\,\partial_{x_i} h(\pmb{x}), \qquad (3.1)$$

the calculation of $\partial_{x_i} Z_n$ requires the calculation of $\partial_\eta Z_n$ and $\partial_h Z_n$, that is, the calculation of the derivatives of $Z_n$ with respect to the boundary points $z = \eta(\pmb{x})$ and $z = h(\pmb{x})$.

In order to keep the presentation as simple as possible without essential loss of generality, we shall proceed by considering the Sturm-Liouville problem (2.3) with $Q = 0$ [9]. On physical grounds, this problem will be referred to as the *water-wave reference waveguide*. For any $\pmb{x}$, the eigenvalues $k_n = k_n(\eta,h) = k_n(H)$, $n \geq 0$ (also called *local wave numbers*), are given as the roots of the transcendental equations

$$\mu_0 - k_0 \tanh(k_0 H) = 0, \qquad \mu_0 + k_n \tan(k_n H) = 0, \quad n \geq 1. \qquad (3.2\text{a,b})$$

The corresponding eigenfunctions $Z_n$, normalized so that $[Z_n]_{z=\eta} = 1$, are given by the formulae

$$Z_0(z;\eta,h) = \frac{\cosh[k_0(z+h)]}{\cosh[k_0(\eta+h)]}, \qquad Z_n(z;\eta,h) = \frac{\cos[k_n(z+h)]}{\cos[k_n(\eta+h)]}. \qquad (3.3\text{a,b})$$

The above eigensystem has been used in several coupled-mode methods in the context of linear water waves [17], [18], [20], that is, when the upper surface of the strip-like domain is flat, $\eta = 0$. The case studied herein corresponds to a general upper surface $z = \eta(\boldsymbol{x})$, as is appropriate for nonlinear water waves; see [25], [27]. Special attention should be paid on the implicit dependence of $Z_n(z;\eta,h)$ on $(\eta,h)$ that comes through the definition of $k_n(H)$, Eqs. (3.2). From now on we shall assume that the boundary functions $\eta$ and $h$ belong to $C^2(\tilde{X}_{\text{part}})$. Then, all order estimates derived below hold true uniformly in $\tilde{X}_{\text{part}}$.

---

[9] No essential loss of generality results from the choice $Q = 0$, in the sense that all the results obtained in this and the following sections remain valid for the case $Q = \text{const.} \neq 0$, the proof being completely similar. New difficulties arise when $Q = Q(\pmb{x},z)$, although the final results, e.g. Theorem 2, are expected to be valid even in this case, under appropriate smoothness assumptions on $Q(\pmb{x},z)$.

**Proposition 1.** *For $n \geq 1$ the first and second derivatives of $k_n(H)$ are given by the following equations, and exhibit the asymptotic behavior indicated therein*:

$$\partial_H k_n = -\frac{k_n(k_n^2 + \mu_0^2)}{-\mu_0 + H(k_n^2 + \mu_0^2)} = O(n), \tag{3.4}$$

$$\partial_H^2 k_n = -2\,\partial_H k_n \left\{ \mu_0 + \frac{\partial_H k_n}{k_n}(H\mu_0 - 1)\left(2 + H\frac{\partial_H k_n}{k_n}\right)\right\} = O(n). \tag{3.5}$$

*Further, by the chain rule we obtain*

$$\partial_{x_i} k_n = \partial_H k_n \,\partial_{x_i} H = O(n), \tag{3.6}$$

$$\partial_{x_i}^2 k_n = \partial_H^2 k_n \,(\partial_{x_i} H)^2 + \partial_H k_n \,\partial_{x_i}^2 H = O(n). \tag{3.7}$$

*Proof*: The first equalities in Eqs. (3.4) and (3.5) are obtained from Eq. (3.2b) by applying the implicit function theorem. The estimates appearing in the rightmost member of Eqs. (3.4), (3.5) are obtained by invoking the estimate $k_n = O(n)$, Eq.(2.5a). Detailed calculations can be found in Appendix A(a). ∎

In order to proceed with the $\boldsymbol{x}$–derivatives of $Z_n = Z_n(z;\eta(\boldsymbol{x}),h(\boldsymbol{x}))$, it is useful to introduce the auxiliary function

$$W_n(z;\eta,h) = W_n(\boldsymbol{x},z) = \frac{\sin[k_n(z+h)]}{\cos[k_n(\eta+h)]} = -\frac{\partial_z Z_n}{k_n}, \tag{3.8}$$

which satisfies the same differential equation as $Z_n$, namely $\partial_z^2 W_n = -k_n^2 W_n$.

**Proposition 2.** *The $\boldsymbol{x}$–derivatives of $Z_n(\boldsymbol{x},z)$, $n \geq 1$, $\boldsymbol{x} \in \tilde{X}_{\text{part}}$, are given by the following equations, and their asymptotic behavior with respect to $n$, as $n \to \infty$, is as indicated therein*:

$$\partial_{x_i} Z_n = -\left(\partial_{x_i} k_n\,(z+h) + k_n\,\partial_{x_i} h\right) W_n - \mu_0 \frac{\partial_{x_i}(k_n H)}{k_n} Z_n = O(n), \tag{3.9}$$

$$\partial_{x_i}^2 Z_n = -(\partial_{x_i} k_n)^2 (z+h)^2 Z_n - 2k_n(\partial_{x_i} k_n)(\partial_{x_i} h)(z+h) Z_n$$
$$+ \left\{-k_n^2(\partial_{x_i} h)^2 + (\partial_{x_i}(k_n H))^2 - \mu_0 \frac{\partial_{x_i}^2(k_n H)}{k_n} + 2\mu_0^2 \left(\frac{\partial_{x_i}(k_n H)}{k_n}\right)^2\right\} Z_n$$
$$- \left\{\partial_{x_i}^2 k_n - 2\mu_0 \frac{\partial_{x_i} k_n}{k_n}\partial_{x_i}(k_n H)\right\}(z+h) W_n$$

$$- \left\{ 2(\partial_{x_i} k_n)(\partial_{x_i} h) + k_n \partial_{x_i}^2 h - 2\mu_0 (\partial_{x_i} h) \partial_{x_i}(k_n H) \right\} W_n = O(n^2). \tag{3.10}$$

*Proof*: The starting point for the derivation of the above equations is the formula (3.1). Derivatives $\partial_\eta Z_n$ and $\partial_h Z_n$ are calculated directly by using Eqs. (3.2b), (3.3b), (3.6) - (3.8). Details are given in Appendix A(b). The estimates in the rightmost member of Eqs. (3.9), (3.10), are derived by applying the estimates given in Proposition 1, and their consequences given in Lemma A1, in Appendix A(b), where all the relevant details are included. ∎

**Remark 3.1.** Formulae for the first derivatives of the eigensystem $\{k_n, Z_n\}$ have also been derived by Porter & Staziker [46], p. 375, in the case of a flat upper surface ($\eta = 0$). Their results are in different form, containing trigonometric functions, but they can be identified with ours, Eqs. (3.4) and (3.9), after some tedious algebraic manipulations. The elimination of trigonometric functions is possible by exploiting the local dispersion relation, Eq. (3.2b).

**3.2.** The derivatives $\partial_{x_i} Z_n$ and $\partial_{x_i}^2 Z_n$, Eqs. (3.9), (3.10), are expressed as linear combinations of $Z_n$, $W_n$, with coefficients which are rational functions of $k_n$. This structure facilitates the derivation of the asymptotic estimates needed for proving the convergence Theorem 2 (section 4).

## 4. Decay rates of derivatives of the modal amplitudes, and termwise differentiability of the series expansion

To establish the termwise differentiability of the coupled-mode series (2.2), corresponding to the water-wave reference waveguide, we shall estimate the decay rates of the differentiated terms. Consider, for example, the terms $T_n^{(\partial x_i)}$, $T_n^{(\partial^2 x_i)}$, obtained after differentiating, with respect to $x_i$, the general term $T_n = \varphi_n Z_n$ of the series (2.2) once or twice, respectively:

$$T_n^{(\partial x_i)}(\boldsymbol{x}, z) = \partial_{x_i} \varphi_n(\boldsymbol{x}) Z_n(\boldsymbol{x}, z) + \varphi_n(\boldsymbol{x}) \partial_{x_i} Z_n(\boldsymbol{x}, z), \tag{4.1}$$

$$T_n^{(\partial^2 x_i)}(\boldsymbol{x}, z) = \partial_{x_i}^2 \varphi_n(\boldsymbol{x}) Z_n(\boldsymbol{x}, z) + \varphi_n(\boldsymbol{x}) \partial_{x_i}^2 Z_n(\boldsymbol{x}, z) + \\ + 2 \partial_{x_i} \varphi_n(\boldsymbol{x}) \partial_{x_i} Z_n(\boldsymbol{x}, z). \tag{4.2}$$

The results obtained in Sec. 2, for $\varphi_n$, and in Sec. 3, for $\partial_{x_i} Z_n$ and $\partial_{x_i}^2 Z_n$, permit us to directly estimate the decay rates of $\varphi_n \partial_{x_i} Z_n$ and $\varphi_n \partial_{x_i}^2 Z_n$, appearing in the right-hand sides of Eqs. (4.1) and (4.2). To estimate the rates of the remaining terms, we need estimates of the derivatives $\partial_{x_i} \varphi_n$ and $\partial_{x_i}^2 \varphi_n$ of the modal amplitudes. Such estimates will be derived in this section by exploiting the representation (2.17) of $\varphi_n$. Eq. (2.17) with $Q_n^* = 1$ is rewritten here in the form

$$\varphi_n(\boldsymbol{x}) = \gamma_n(\boldsymbol{x}) k_n^{-4}(\boldsymbol{x}) \lambda_n(\boldsymbol{x}), \tag{4.3}$$

where $\gamma_n(\boldsymbol{x}) = \|Z_n(\boldsymbol{x}, \cdot)\|_{L^2}^{-2}$, and

$$\lambda_n(\boldsymbol{x}) = \left[\partial_z^2 \Phi^* \mu_0 - \partial_z^3 \Phi^*\right]_{z=\eta} - \left[\partial_z^3 \Phi^*\right]_{z=-h} Z_n(-h) + \int_{-h}^{\eta} \partial_z^4 \Phi^* Z_n \, dz, \qquad (4.4)$$

with $\Phi^*$ given by Eq. (2.9c). Clearly, we need estimates of the quantities $\gamma_n(\boldsymbol{x})$, $k_n^{-4}(\boldsymbol{x})$, $\lambda_n(\boldsymbol{x})$ and $Z_n(-h)$, and their $\boldsymbol{x}$–derivatives.

**Lemma 1.** *Assume that $\eta$ and $h$ belong to $C^2(\tilde{X}_{\text{part}})$. Then, the following estimates hold true uniformly in $\tilde{X}_{\text{part}}$:*

$$k_n^{-4}(\boldsymbol{x}), \ \partial_{x_i} k_n^{-4}(\boldsymbol{x}), \ \partial_{x_i}^2 k_n^{-4}(\boldsymbol{x}) \text{ are all of order } O(n^{-4}). \qquad (4.5)$$

$$\gamma_n(\boldsymbol{x}), \ \partial_{x_i} \gamma_n(\boldsymbol{x}), \ \partial_{x_i}^2 \gamma_n(\boldsymbol{x}) \text{ are all of order } O(1). \qquad (4.6)$$

$$Z_n(-h) = O(1), \ \partial_{x_i} Z_n(-h) = O(n^{-2}), \ \partial_{x_i}^2 Z_n(-h) = O(n^{-2}). \qquad (4.7)$$

*Proof*: See Appendix B(a).

The following lemma provides estimates for integrals involved in the definition of $\lambda_n(\boldsymbol{x})$ and its derivatives $\partial_{x_i} \lambda_n(\boldsymbol{x})$, $\partial_{x_i}^2 \lambda_n(\boldsymbol{x})$.

**Lemma 2:** *Let $F$ be a function defined on $\bar{D}_h^\eta(X)$, and assume that $F(\boldsymbol{x}, \cdot)$ satisfies the local $VSC(k=2; \boldsymbol{x})$. Then, locally at $\boldsymbol{x}$, we have*

$$\text{(i)} \int_{-h}^{\eta} F Z_n \, dz = O(n^{-2}), \qquad \text{(ii)} \int_{-h}^{\eta} F W_n \, dz = O(n^{-1}). \qquad (4.8\text{a,b})$$

*If, in addition, $h, \eta \in C^2(\tilde{X}_{\text{part}})$, then*

$$\text{(iii)} \int_{-h}^{\eta} F \partial_{x_i} Z_n \, dz = O(1), \qquad \text{(iv)} \int_{-h}^{\eta} F \partial_{x_i}^2 Z_n \, dz = O(1). \qquad (4.9\text{a,b})$$

*Further, if $F(\boldsymbol{x}, \cdot)$ satisfies the uniform $VSC(k=2; \tilde{X}_{\text{part}})$, then, the estimates* (4.8) *and* (4.9) *are uniformly valid in $\tilde{X}_{\text{part}}$.*

*Proof:* Assertions **(i)** and **(ii)** follow easily by using equations $Z_n = -k_n^{-2} \partial_z^2 Z_n$ and $W_n = -k_n^{-2} \partial_z^2 W_n$, performing two integrations by parts (permitted by the local $VSC(k=2; \boldsymbol{x})$), and estimating the arising terms. The main idea for proving assertions **(iii)** and **(iv)**, is to express $\partial_{x_i} Z_n$ and $\partial_{x_i}^2 Z_n$ as linear combinations of $Z_n$ and $W_n$, using Eqs. (3.9) and (3.10), and then estimate the arising terms by exploiting assertions **(i)** and **(ii)**, Eq. (2.5a), and Proposition 1 together with its consequence given in Lemma A1. The details of the above procedure are presented in Appendix B(b). ∎

**Lemma 3.** *If the field $\Phi(\boldsymbol{x}, \cdot)$ and its derivatives $\partial_{x_i} \Phi(\boldsymbol{x}, \cdot)$ and $\partial_{x_i}^2 \Phi(\boldsymbol{x}, \cdot)$ satisfy the uniform $VSC(k=6; \tilde{X}_{\text{part}})$, then*

$$\lambda_n(\boldsymbol{x}), \ \partial_{x_i}\lambda_n(\boldsymbol{x}), \ \partial^2_{x_i}\lambda_n(\boldsymbol{x}) \ \text{are all of order} \ O(1). \tag{4.10}$$

*Proof*: The result $\lambda_n(\boldsymbol{x}) = O(1)$ follows from Eq. (4.4). Indeed, from the assumption $\Phi \in H^6(-h,\eta)$ and Eq. (2.9c), we conclude that $\Phi^* \in H^6(-h,\eta)$, and thus all $\Phi^*$−boundary values appearing in the right-hand side of Eq. (4.4) are well-defined and $O(1)$. Since $Z_n(-h) = O(1)$, according to Eq. (4.7), all boundary terms are of order $O(1)$. The integral term is found to be $O(n^{-2})$ by using Eq. (4.8a) with $F = \partial_z^4 \Phi^* \in H^2(-h,\eta)$.

To prove $\partial_{x_i}\lambda_n = O(1)$ and $\partial^2_{x_i}\lambda_n = O(1)$, we first calculate $\partial_{x_i}\lambda_n$ and $\partial^2_{x_i}\lambda_n$ by differentiating Eq. (4.4); see Eqs. (B9), (B10) in Appendix B(c). Smoothness assumptions and Eq. (4.7) ensure that all boundary terms are of order $O(1)$. The various integral terms are estimated by applying Lemma 2 with appropriate choices of $F$. The details are given in Appendix B(c). ∎

The estimates (4.5), (4.6) and (4.10), in conjunction with Eq. (4.3), lead to the following

**Theorem 2.** *Assume that* $h, \eta \in C^2(\tilde{X}_{\text{part}})$, *and* $\Phi(\boldsymbol{x}, \cdot)$ *is as in Lemma 3. Then, the sequences* $\partial_{x_i}\varphi_n$ *and* $\partial^2_{x_i}\varphi_n$ *satisfy the following estimates with respect to* $n$, *uniformly on* $\tilde{X}_{\text{part}}$:

i) $\left\| \partial_{x_i} \varphi_n \right\|_{\infty, \tilde{X}_{\text{part}}} = \sup\{|\partial_{x_i}\varphi_n(\boldsymbol{x})|, \ \boldsymbol{x} \in \tilde{X}_{\text{part}}\} = O(n^{-4}),$ (4.11)

ii) $\left\| \partial^2_{x_i} \varphi_n \right\|_{\infty, \tilde{X}_{\text{part}}} = \sup\{|\partial^2_{x_i}\varphi_n(\boldsymbol{x},t)|, \ \boldsymbol{x} \in \tilde{X}_{\text{part}}\} = O(n^{-4}).$ (4.12)

*If the smoothness assumptions hold on* $X$, *then* $\tilde{X}_{\text{part}}$ *can be replaced by* $X$ *in all statements of the theorem.*

**Corollary 2.** *Under the assumptions stated in Theorem 2,* $\left\|\partial_\alpha(\varphi_n Z_n)\right\|_\infty = O(n^{-3})$ *and* $\left\|\partial^2_\alpha(\varphi_n Z_n)\right\|_\infty = O(n^{-2})$ *for* $\alpha \in \{x_1, x_2, z\}$. *Thus, the series expansion* (2.2) *can be differentiated term-by-term twice with respect to all variables* $x_1, x_2$, *and* $z$. *The termwise differentiated series are uniformly convergent in* $\bar{D}_h^\eta(\tilde{X}_{\text{part}})$ *or in* $\bar{D}_h^\eta(X)$, *if the smoothness assumptions hold true uniformly in* $X$.

## 5. Application to an elliptic boundary-value problem

In this section, we apply the exact semi-separation of variables (2.2) to the solution of the following elliptic BVP in a strip-like domain $D_h^\eta(X)$:

$$\Delta\Phi = 0, \quad \text{in} \quad D_h^\eta(X), \tag{5.1a}$$

$$[\Phi]_{z=\eta} = \psi, \qquad [\partial_n \Phi]_{z=-h} = 0, \tag{5.1b,c}$$

and periodic lateral-boundary conditions. As it is shown in Appendix C(a), using the exact expansion $\Phi = \sum_{n=-2}^{\infty} \varphi_n(\boldsymbol{x}) \, Z_n(z; \eta(\boldsymbol{x}), h(\boldsymbol{x}))$ in conjunction with a variational formulation of problem (5.1), we obtain the following equivalent consistent[10] coupled-mode system (CCMS):

$$\sum_{n=-2}^{\infty} \left( A_{mn} \nabla_x^2 + \boldsymbol{B}_{mn} \cdot \nabla_x + C_{mn} \right) \varphi_n = 0, \quad m \geq -2, \quad \boldsymbol{x} \in X, \tag{5.2a}$$

$$\sum_{n=-2}^{\infty} \varphi_n = \psi, \quad \text{on } X. \tag{5.2b}$$

The $\boldsymbol{x}$-dependent coefficients $A_{mn}$, $\boldsymbol{B}_{mn}$, $C_{mn}$ are given by the formulae

$$A_{mn} = \int_{-h}^{\eta} Z_n Z_m \, dz, \tag{5.3a}$$

$$\boldsymbol{B}_{mn} = 2 \int_{-h}^{\eta} (\nabla_x Z_n) Z_m \, dz + (\nabla_x h) \left[ Z_m Z_n \right]_{z=-h}, \tag{5.3b}$$

$$C_{mn} = \int_{-h}^{\eta} (\nabla_x^2 Z_n + \partial_z^2 Z_n) Z_m \, dz - \boldsymbol{N}_h \cdot \left[ (\nabla_x Z_n, \partial_z Z_n) Z_m \right]_{z=-h}, \tag{5.3c}$$

where $\boldsymbol{N}_h = (-\nabla_x h, -1)$ is the outward (with respect to $D_h^\eta$) normal (but not normalized to have unit length) vector on the lower boundary $\Gamma_h$. Since the series representation (2.2) is not subjected to any restriction concerning the behavior of $\Phi$ (apart from smoothness requirements), the above reformulation remains valid for any type of lateral conditions; see Appendix C(a,c). For computing a numerical solution, we choose the total number of modes $N_{tot} = M + 3$ [11] to be retained in the expansion (2.2), truncating the modal sequence to $\{\varphi_n(\boldsymbol{x})\}_{n=-2}^{M}$. These $\varphi_n(\boldsymbol{x})$ are calculated by using the first $N_{tot} - 1$ partial differential equations (PDEs) (5.2a), together with the linear algebraic constraint (5.2b) [12]. The obtained system is discretized via a finite-difference (FD) method (here we use fourth-order central FD), leading to a square, sparse, linear algebraic system, which is solved numerically.

In order to illustrate the convergence and the accuracy of the present approach, we consider the one dimensional version of problem (5.1) ($\nabla_x \equiv (\partial_x, 0)$) in the laterally-periodic case, $X = [0, 2\pi]$, $\Phi(0, z) = \Phi(2\pi, z)$, $\partial_x \Phi(0, z) = \partial_x \Phi(2\pi, z)$. Especially, we investigate numerically two examples, introduced by Nicholls and Reitich [47], for which analytical solutions are known. Consider the function $\Phi_\kappa(x, z) = \cosh(\kappa(z + h_0)) \cos(\kappa x)$ in the domain $0 \leq x \leq 2\pi$, $-h_0 \leq z \leq \eta(x)$. $\Phi_\kappa(x, z)$ satisfies Eqs. (5.1a) and (5.1c) at $z = -h_0$, while the values $\left[ \Phi \right]_{z=\eta} = \psi$ can be calculated directly for any choice of upper surface $z = \eta(x)$

---

[10] The term *consistent* refers to the fact that the infinite series appearing in Eqs. (5.2) are convergent.

[11] $M$ denotes the order of the last mode retained in the truncated version of expansion (2.2).

[12] Other truncation strategies are possible and will be discussed elsewhere.

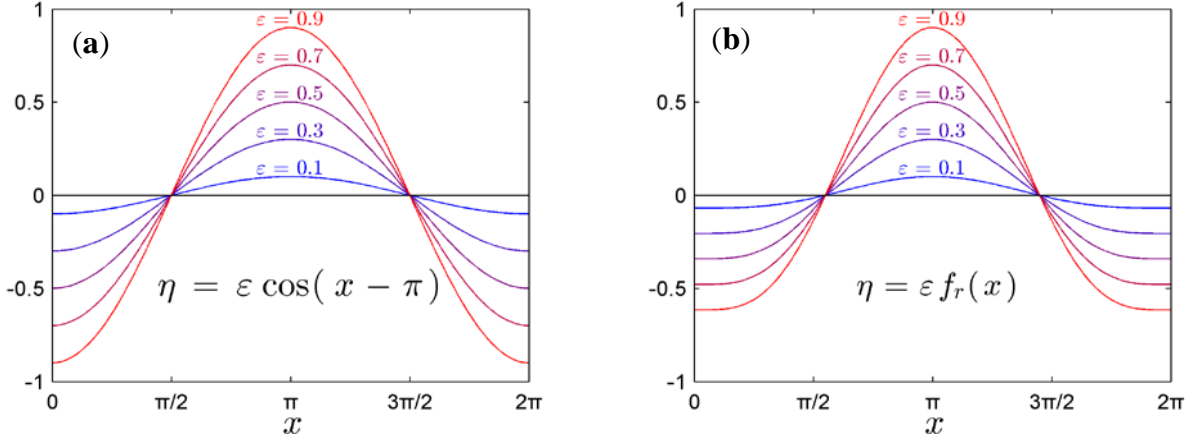

**Figure 2.** The periodic domains used in the computations. (**a**) Smooth case, (**b**) Rough case.

Two families of upper-surface profiles are introduced: $\eta(x) = \varepsilon \cos(\kappa(x+\gamma))$, where $\gamma$ is an indifferent constant, and $\eta(x) = \varepsilon f_r(r)$, $f_r(x) = A x^4 (2\pi - x)^4 + B$, where the constants $A$ and $B$ are chosen so that $\int_0^{2\pi} f_r(x)\,dx = 0$ and $\max_{x \in [0, 2\pi]} \{f_r(x)\} = f_r(\pi) = 1$. Following [47], we call the surfaces of the first family *smooth* ($\eta \in C^\infty(\mathbb{R})$), and the ones of the second family *rough* ($\eta \in C^4(\mathbb{R})$). Computations presented below correspond to $\kappa = 1$, $h_0 = 1$ and various values of the deformation parameter $\varepsilon$, ranging from $\varepsilon = 0.1$ (mildly deformed domains) to $\varepsilon = 0.9$ (strongly deformed domains); see **Figure 2**. For the construction of the vertical basis functions $Z_n$, we take $\mu_0 = \kappa \tanh(\kappa h_0)$. The $x-$discretization is performed by using a uniform grid of $256$ points.

**(a) Convergence of the numerical solution to the exact solution**

As a first test, system (5.2) is truncated at a large number of modes $N_{tot} = 70$, and solved numerically, in order to investigate the asymptotic behaviour of numerically computed $\varphi_n$, $\partial_x \varphi_n$, $\partial_x^2 \varphi_n$, for large values of $n$. Results on the uniform $C^2$−norm, $\|\varphi_n\|_{C^2(X)} = \|\varphi_n\|_\infty + h_0 \|\partial_x \varphi_n\|_\infty + h_0^2 \|\partial_x^2 \varphi_n\|_\infty$, for both the smooth and rough cases, and three values of the deformation parameter, $\varepsilon = 0.1, 0.5, 0.9$, are shown in **Figure 3**. In both cases, the theoretical rate of decay, $\|\varphi_n\|_{C^2(X)} = O(n^{-4})$ (proved in Theorems 1 and 2), is accurately reproduced by the numerical solution. Most importantly, it is observed that the decay of the first few $\varphi_n$ is in fact exponential. This phenomenon is an inherent, generic feature of the present method, appearing consistently in all geometric configurations. Although a theoretical proof is lacking, it is reasonable to assume that this "superconvergence" plays an essential role in the efficiency of the present method.

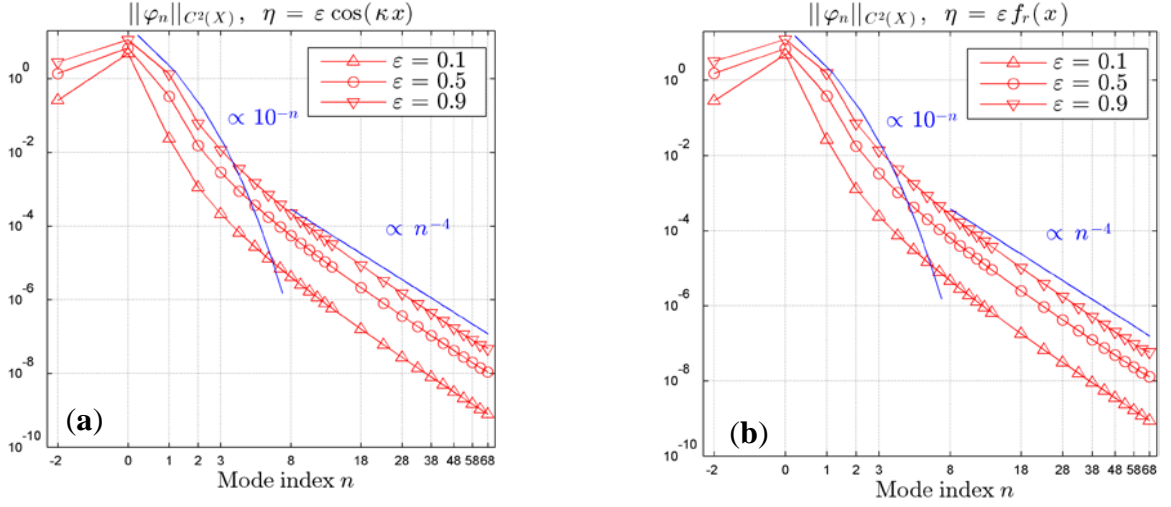

**Figure 3.** Decay rate of numerically computed $\|\varphi_n\|_{C^2(X)}$. (**a**) Smooth case, (**b**) Rough case.

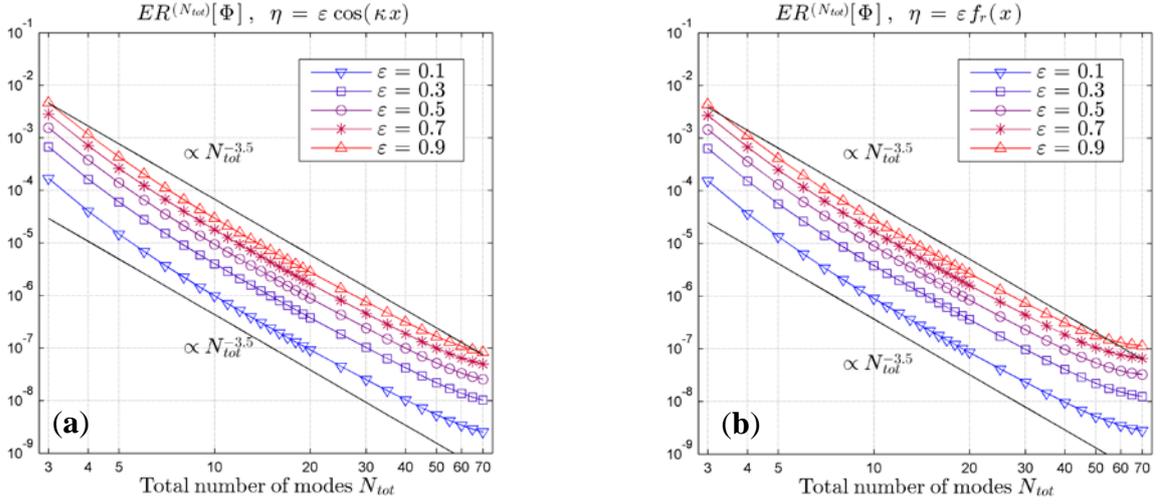

**Figure 4.** $L^2$ − relative error of the solution vs. the total number of modes $N_{tot}$ used, for $\varepsilon = 0.1(0.2)0.9$. (**a**) Smooth case, (**b**) Rough case.

The next and more important test concerns the convergence rate and the accuracy of the numerical solution of the whole field $\Phi$ with respect to the order of truncation $N_{tot}$. Numerical results on the $L^2$ relative error $ER^{(N_{tot})}[\Phi] := \left\| \Phi_\kappa - \sum_{n=-2}^{M} \varphi_n Z_n \right\|_{L^2(D_h^\eta)} / \|\Phi_\kappa\|_{L^2(D_h^\eta)}$ (recall that $\Phi_\kappa$ is the exact solution of the studied problem) have been obtained for $\varepsilon = 0.1(0.2)0.9$, by using $N_{tot} = 3,...,70$, and are shown in **Figure 4**. In all cases, the rate of decay of the error is at worst proportional to $N_{tot}^{-3.5}$, for all values of the deformation parameter $\varepsilon$. This rate of convergence of the solution field is reported in several other circumstances (PDEs, boundary conditions) where enhanced eigenfunctions expansions are used [36], [8], [37], [9].

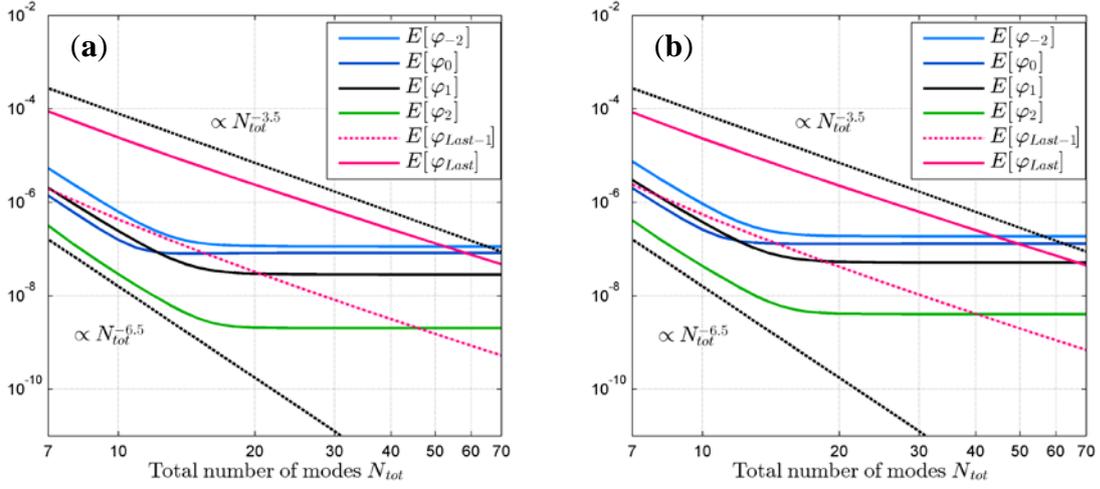

**Figure 5.** $L^2$ – abloulute error of the modal amplitudes vs. the total number of modes $N_{tot}$ for $\varepsilon = 0.5$. (**a**) Smooth case, (**b**) Rough case.

In order to gain more insight on this rate of convergence, we also consider the modal amplitudes $\{\varphi_n^{(N_{tot})}\}_{n=-2}^M$, as calculated by solving numerically Eqs. (5.2) using $N_{tot}$ modes, and compare them with their exact values, obtained by applying Eqs. (2.9) (Theorem 1) to the function $\Phi_\kappa$; see Eqs. (D1), (D2), (D6), (D7) in Appendix D. Since $\varphi_n$ decays rapidly to zero, we consider the absolute $L^2$ – error $E[\varphi_n] \equiv E^{(N_{tot})}[\varphi_n] = \|\varphi_n - \varphi_n^{(N_{tot})}\|_{L^2(X)}$. Results for the first few modal amplitudes, indexed by $n = -2, 0, 1, 2$, and for the two last ones, $\varphi_{Last} \equiv \varphi_M = \varphi_{N_{tot}-2}$ and $\varphi_{Last-1} \equiv \varphi_{M-1} = \varphi_{N_{tot}-3}$, are shown in **Figure 5.** Surprisingly, it is observed that the error of the mode $\varphi_{-2}$ and the first few $\varphi_n$, $n = 0, 1, 2$, decays at a much faster([13]) rate, namely as $N_{tot}^{-6.5}$. This decay rate persists for most of the modal amplitudes, up to the last few ones (three to four), for which the error decay gradually decelerates to $N_{tot}^{-3.5}$. The later rate dominates the rate of decay of the error for the whole solution $\Phi$, as shown in Figure 4. This behaviour can be attributed to the truncation strategy applied to Eqs. (5.2), and deserves further investigation. It can be anticipated that, by optimizing the truncation strategy, the rate of decay of the error for the whole solution would be improved up to $N_{tot}^{-6.5}$. In the next subsection, we shall study another interesting problem, exploiting the superconvergence of the modal amplitude $\varphi_{-2}$.

---

([13]) In comparison with the decay rate of the error of the whole solution, Figure 4.

**(b) Computation of the Dirichlet-to-Neumann operator corresponding to Eqs. (5.1)**

In various applications it is necessary to calculate the Dirichlet-to-Neumann (DtN) operator associated with the BVP (5.1). Given $\eta$, $h$ and $\psi$, the DtN operator is defined by the relation

$$\mathsf{G}[\eta,h]\psi = [\partial_n \Phi]_{z=\eta} \left(1 + |\nabla_x \eta|^2\right)^{1/2}, \tag{5.4}$$

where $\Phi$ is the solution of Eqs. (5.1). This operator plays a fundamental role in the nonlinear water-wave problem, where $\eta$ is the free-surface elevation, $h$ is the seabed profile, and $\Phi$ is the velocity potential in the fluid domain $D_h^\eta(X)$, having $\psi$ as trace on the free surface. In fact, the Hamiltonian evolution equations of this problem are formulated in terms of the DtN operator as follows (see, e.g., [48], [41], [49])

$$\partial_t \eta = \mathsf{G}[\eta,h]\psi,$$

$$\partial_t \psi = -g\eta - \frac{1}{2}(\nabla_x \psi)^2 + \frac{(\mathsf{G}[\eta,h]\psi + \nabla_x \eta \cdot \nabla_x \psi)^2}{2(|\nabla_x \eta|^2 + 1)}.$$

Exploiting the exact semi-separation of variables (2.2), it can be shown that the DtN operator (5.4) is expressed by the following simple formulae, in terms of the data $\eta$, $\psi$ and the single modal amplitude $\varphi_{-2}$:

$$\mathsf{G}[\eta,h]\psi = -(\nabla_x \eta) \cdot (\nabla_x \psi) + (|\nabla_x \eta|^2 + 1)\left(h_0^{-1} \varphi_{-2} + \mu_0 \psi\right). \tag{5.5}$$

See Appendix C(b) for the derivation of Eq. (5.5). Clearly, the rapid convergence of $\varphi_{-2}^{(N_{tot})}$ is inherited to the approximation of the DtN operator furnished by Eq. (5.5). This is of great practical importance for the numerical treatment of nonlinear water waves, where highly accurate approximations of $\mathsf{G}[\eta,h]\psi$ are needed. Numerical results on $\mathsf{G}^{(N_{tot})}[\eta,h]\psi$ for the cases studied above, with $\varepsilon = 0.1(0.2)0.9$, are given and discussed in Appendix E(a). As expected, the $L^2$ − relative error of $\mathsf{G}[\eta,h]\psi$, decays at a rate proportional to $N_{tot}^{-6.5}$, independently of the value of the deformation parameter $\varepsilon$ and the smoothness of upper surface ($C^\infty$ or $C^4$). Moreover, in order to illustrate the efficiency of the present method for the case where both upper and lower boundaries are non-planar, we have computed the DtN operator for a strip-like domain with strongly "undulated" boundaries, ending to a flat part with different left and right depths; see Appendix E(b). Comparison of numerical solutions obtained by using the present method and a direct numerical method (Boundary Element Method) shows that highly accurate computations of the DtN operator are achieved by the exact semi-separation of variables, at about $10^2$ times faster than by the direct numerical method.

**6. Discussion and conclusions**

In this paper we have established the validity of the exact semi-separation of variables, Eq. (2.2), throughout any strip-like domain with smooth, non-planar boundaries, and any sufficiently smooth field $\Phi$. The set of transverse (vertical) basis functions $\{Z_n, n \geq -2\}$ consists of a local $L^2$ − basis $\{Z_n, n \geq 0\}$ augmented by two additional modes, making it an $H^2$–basis. The convergence properties of the series expansion are not affected by the shape (slope, mean level, amplitude variation) of the boundaries, assuming they remain smooth. No boundary conditions are imposed *a priori* on the function-to-be-expanded. The traces of $\Phi$ and its first derivatives $\partial_\alpha \Phi$, $\alpha \in \{x_1, x_2, z\}$, on the boundary $\partial D_h^\eta(X)$ may be any continuous functions. The expansion (2.2) may be used as the starting point for formulating exact CMTs (called herein consistent CMTs, CCMTs) for any second-order operator acting on functions defined on $D_h^\eta(X)$, with general (conventional or unconventional) first-order boundary conditions. An interesting example is

$$\nabla \cdot (\rho(\boldsymbol{x}, z) \nabla \Phi) + \sigma^2(\boldsymbol{x}, z) \Phi = 0, \quad (\boldsymbol{x}, z) \in D_h^\eta(X), \tag{6.1a}$$

$$\left[\partial_n \Phi - \mu_0(\boldsymbol{x}) \Phi\right]_{z=\eta(\boldsymbol{x})} = \mu_1(\boldsymbol{x}), \tag{6.1b}$$

$$\left[\partial_n \Phi - \nu_0(\boldsymbol{x}) \Phi\right]_{z=-h(\boldsymbol{x})} = \nu_1(\boldsymbol{x}), \tag{6.1c}$$

where $\rho(\boldsymbol{x}, z)$, $\sigma^2(\boldsymbol{x}, z)$, $\mu_j(\boldsymbol{x})$, $\nu_j(\boldsymbol{x})$, $j = 0, 1$, are sufficiently smooth functions. In these applications the $L^2$–sub-basis $\{Z_n, n \geq 0\}$ may be constructed by the eigenfunctions of the reference waveguide, making the CCMTs exact generalizations of the usual normal-mode solutions to the case of irregular strip-like domains. In this way, a great deal of the physical structure carried by normal modes is retained and, at the same time, the modifications of the field due to the boundary deformations are correctly taken into account.

It is important to stress here that the semi-separation of variables, Eq. (2.2), and the implied CCMT *is not a perturbative approach*. The truncation used in implementing numerical solution schemes is realized by ignoring terms (of the expansion (2.2)) which are consistently small, of orders $O(n^{-4})$, independently of how much the boundaries differ from the planar ones. Thus, even the finite-dimensional approximations obtained by a CCMT *are not perturbative*, as regards the boundary shape. This feature distinguishes CCMT from various other methods which are based on perturbation expansions and, thus, become inefficient and eventually divergent for high values of the slope or amplitude of the boundary deformation. Besides, the present approach realizes a dimension reduction of the problem, providing results having accuracy and stability comparable with those of direct numerical methods, being much lighter computationally. The results presented in this paper, lead to the conclusion that the exact semi-separation of variables, proposed herein, comprises a stable and efficient alternative to pertubative and direct numerical methods for the solution of boundary value problems in general strip-like domains. The proposed method retains the efficiency of the pertubative methods and the accuracy of the direct ones, exhibiting a striking superconvergence as regards the computation of the corresponding DtN operators.

**Data accessibility.** This paper is theoretical and includes in the text all necessary data and equations to support the results presented. No experimental data has been used in this work.

**Authors' contributions.** GAA coordinated the study, and supervised the paper writing; both authors contributed to the development of the theoretical part; Ch.EP conducted the numerical calculations; both authors equally contributed to the paper writing and gave final approval for publication.

**Competing interests.** We have no competing interests.

**Funding.** This work is a result of unsponsored Academic research.

**Acknowledgments:** Ch.E. Papoutsellis was supported by a PhD scholarship from NTUA.

## Appendix A. Auxiliary results for Section 3

All calculations presented in this and the following appendices pertain to the constant-coefficient, water-wave reference waveguide, which is obtained from Eqs. (2.3) with $Q = 0$. For easy reference we write down the formulation of the Sturm-Liouville problem corresponding to this waveguide:

$$\partial_z^2 Z_n(z) + k_n^2 Z_n(z) = 0, \qquad -h(\boldsymbol{x}) < z < \eta(\boldsymbol{x}), \tag{A1a}$$

$$B^\eta Z_n \equiv [(\partial_z - \mu_0) Z_n]_{z=\eta} = 0, \tag{A1b}$$

$$B_h Z_n \equiv [\partial_z Z_n]_{z=-h} = 0. \tag{A1c}$$

We shall now give detailed proofs of Propositions 1, and 2, stated in section 3 of the main paper.

### (a) Detailed Proof of Proposition 1

Setting

$$G_n(H, k_n) = \mu_0 + k_n \tan(k_n H), \tag{A2}$$

the local dispersion relation, implicitly defining the function $k_n = k_n(H)$, is written in the form $G_n(H, k_n) = 0$. The implicit function theorem implies that

$$\partial_H k_n = -\partial_H G_n / \partial_{k_n} G_n, \tag{A3}$$

$$\partial_H^2 k_n = -\frac{\partial_H^2 G_n + 2(\partial_{H k_n} G_n) \partial_H k_n + \partial_{k_n}^2 G_n (\partial_H k_n)^2}{\partial_{k_n} G_n}. \tag{A4}$$

Straightforward calculations yield

$$\partial_H G_n = k_n^2 \sec^2(k_n H), \tag{A5a}$$

$$\partial_{k_n} G_n = \tan(k_n H) + k_n H \sec^2(k_n H), \tag{A5b}$$

which, in conjunction with Eq. (A3), leads to

$$\partial_H k_n = -\frac{k_n^2 \sec^2(k_n H)}{\tan(k_n H) + k_n H \sec^2(k_n H)}. \tag{A6}$$

Eliminating the trigonometric functions in the right-hand side of the above equation by using Eq. (3.2b) and the identity $1 + \tan^2 a = \sec^2 a$, Eq. (A6) results in

$$\partial_H k_n = -\frac{k_n(k_n^2 + \mu_0^2)}{-\mu_0 + H(k_n^2 + \mu_0^2)}, \tag{A7}$$

which proves the first equality in Eq. (3.4). The second equality in Eq. (3.4), that is the estimate $\partial_H k_n = O(n)$, follows easily from the estimate $k_n = O(n)$, Eq. (2.5a), and the following limiting relation, obtained by using Eq. (A7),

$$\frac{\partial_H k_n}{k_n} = -\frac{1 + \mu_0^2 / k_n^2}{H - (\mu_0 - H\mu_0^2)/k_n^2} \to -\frac{1}{H}, \quad \text{as } n \to \infty. \tag{A8}$$

The proof of Eq. (3.4) is now completed.

Eq. (3.5) is proved by using Eq. (A4). Differentiating twice $G_n(H, k_n)$ and using Eq. (3.2b), and the trigonometric identity $1 + \tan^2 a = \sec^2 a$, we get

$$\partial_H^2 G_n = 2k_n^3 \sec^2(k_n H) \tan(k_n H) = -2\mu_0 \left(k_n^2 + \mu_0^2\right), \tag{A9a}$$

$$\partial_{Hk_n} G_n = 2k_n \sec^2(k_n H) + 2k_n^2 H \sec^2(k_n H) \tan(k_n H) =$$
$$= 2\frac{k_n^2 + \mu_0^2}{k_n}\left(1 - H\mu_0\right), \tag{A9b}$$

$$\partial_{k_n}^2 G_n = 2H \sec^2(k_n H) + 2k_n H^2 \sec^2(k_n H) \tan(k_n H) =$$
$$= 2H\frac{k_n^2 + \mu_0^2}{k_n^2}\left(1 - H\mu_0\right). \tag{A9c}$$

Substituting the derivatives of $G_n(H, k_n)$ from Eqs. (A5b) and (A9) into Eq. (A4), and using Eq. (A7), we find

$$\partial_H^2 k_n = -2\,\partial_H k_n \left\{ \mu_0 + \frac{\partial_H k_n}{k_n}\left(H\mu_0 - 1\right)\left[2 + H\frac{\partial_H k_n}{k_n}\right]\right\}, \tag{A10}$$

which proves the first equality of Eq. (3.5). The second equality in Eq. (3.5), that is the estimate $\partial_H^2 k_n = O(n)$, follows directly from Eqs. (A10) and (A8), in conjunction with the already proved result $\partial_H k_n = O(n)$. This completes the proof of Proposition 1.

### (b) Detailed proof of Proposition 2

- *Calculation of the first derivative $\partial_{x_i} Z_n$*

The eigenfunctions $Z_n$ are dependent on $\mathbf{x}$ through their dependence on the boundary points (functions) $\eta = \eta(\mathbf{x})$ and $h = h(\mathbf{x})$. Thus, we have

$$\partial_{x_i} Z_n(\mathbf{x}, z) = \partial_\eta Z_n\, \partial_{x_i} \eta(\mathbf{x}) + \partial_h Z_n\, \partial_{x_i} h(\mathbf{x}). \tag{A11}$$

$\partial_\eta Z_n$ and $\partial_h Z_n$ are the first Fréchet (shape) derivatives of the eigenfunctions $Z_n$ with respect to the end points $\eta$ and $h$ of the interval in which they are defined. These derivatives will be now calculated by differentiating the explicit form of $Z_n$ given by Eq. (3.3b). Since $k_n = k_n(\eta, h)$ (defined implicitly through Eq. (3.2b)), it is expedient to use here the notation $Z_n = Z_n(z; \eta, h; k_n(\eta, h))$, in order to keep track of all dependences either explicit or implicit. Thus, recalling that $\eta + h = H$, we rewrite Eq. (3.3b) in the form

$$Z_n = Z_n(z;\eta,h;k_n(\eta,h)) = \cos\left[k_n(\eta,h)(z+h)\right]\cdot\sec\left[k_n(\eta,h)H\right]. \quad (A12)$$

Now, by direct differentiation and application of the chain rule, and using the definition of $W_n$, Eq. (3.8), Eq. (A12) yields

$$\begin{aligned}\partial_\eta Z_n &= -\sin\left[k_n(\eta,h)(z+h)\right]\cdot\partial_\eta k_n\cdot(z+h)\cdot\sec\left[k_n(\eta,h)H\right] + \\ &+ \cos\left[k_n(\eta,h)(z+h)\right]\cdot\tan\left[k_n(\eta,h)H\right]\cdot\sec\left[k_n(\eta,h)H\right]\cdot\left(\partial_\eta k_n\cdot H + k_n\right) = \\ &= -\partial_\eta k_n\cdot(z+h)\cdot W_n + \tan\left[k_n(\eta,h)H\right]\cdot\left(\partial_\eta k_n\cdot H + k_n\right)\cdot Z_n, \quad (A13)\end{aligned}$$

and

$$\begin{aligned}\partial_h Z_n &= -\sin\left[k_n(\eta,h)(z+h)\right]\cdot\left(\partial_h k_n\cdot(z+h) + k_n\right)\cdot\sec\left[k_n(\eta,h)H\right] + \\ &+ \cos\left[k_n(\eta,h)(z+h)\right]\cdot\tan\left[k_n(\eta,h)H\right]\cdot\sec\left[k_n(\eta,h)H\right]\cdot\left(\partial_h k_n\cdot H + k_n\right) = \\ &= -\left(\partial_h k_n\cdot(z+h) + k_n\right)\cdot W_n + \tan\left[k_n(\eta,h)H\right]\cdot\left(\partial_h k_n\cdot H + k_n\right)\cdot Z_n. \quad (A14)\end{aligned}$$

Recalling that $k_n = k_n(\eta,h) = k_n(\eta+h) = k_n(H)$ (see Eq. (3.2b)), we see that

$$\partial_\eta k_n = \partial_h k_n = \partial_H k_n, \quad (A15)$$

which, upon substitution in Eqs. (A13) and (A14), leads to the following equations

$$\partial_\eta Z_n = -\partial_H k_n (z+h) W_n + \tan(k_n H)\left(\partial_H k_n H + k_n\right) Z_n, \quad (A16)$$

$$\partial_h Z_n = -\left(k_n + \partial_H k_n (z+h)\right) W_n + \tan(k_n H)\left(\partial_H k_n H + k_n\right) Z_n. \quad (A17)$$

Substituting $\partial_\eta Z_n$ and $\partial_h Z_n$ from Eqs. (A16) and (A17) into Eq. (A11), and noting that

$$\partial_\eta k_n\cdot\partial_{x_i}\eta + \partial_h k_n\cdot\partial_{x_i}h = \partial_{x_i}k_n \quad \text{and} \quad \partial_{x_i}k_n\cdot H + k_n\cdot\partial_{x_i}H = \partial_{x_i}(k_n H),$$

we obtain

$$\partial_{x_i}Z_n = -\left(\partial_{x_i}k_n\cdot(z+h) + k_n\cdot\partial_{x_i}h\right)\cdot W_n + \tan(k_n H)\cdot\partial_{x_i}(k_n H)\cdot Z_n. \quad (A18)$$

Substituting $\tan(k_n H)$ in Eq. (A18) by $-\mu_0/k_n$ (from Eq. (3.2b)) we directly get the first equality in Eq. (3.9).

- *Calculation of the second derivative $\partial_{x_i}^2 Z_n$*

The calculation of the second derivative, $\partial_{x_i}^2 Z_n$ is more involved. A general formula is obtained by differentiating both members of Eq. (A11) with respect to $x_i$, which reads as follows

$$\begin{aligned}\partial_{x_i}^2 Z_n &= \partial_\eta^2 Z_n (\partial_{x_i}\eta)^2 + \partial_h^2 Z_n (\partial_{x_i}h)^2 + 2\partial_{\eta h}Z_n (\partial_{x_i}\eta)(\partial_{x_i}h) \\ &\quad + \partial_\eta Z_n\,\partial_{x_i}^2\eta + \partial_h Z_n\,\partial_{x_i}^2 h.\end{aligned} \quad (A19)$$

The first Fréchet derivatives, $\partial_\eta Z_n$ and $\partial_h Z_n$, of the eigenfunctions $Z_n$ have been already calculated, Eqs. (A16) and (A17). The second Fréchet derivatives, $\partial_\eta^2 Z_n$, $\partial_h^2 Z_n$ and $\partial_{\eta h} Z_n$, can be calculated by differentiating both members of Eqs. (A16) and (A17) with respect to $\eta$ and $h$. However, it is more convenient (and shorter) to proceed differently, differentiating directly with respect to $x_i$ both members of Eq. (A18). This approach leads to

$$\partial_{x_i}^2 Z_n = -\left((z+h)\partial_{x_i}^2 k_n + 2(\partial_{x_i} h)(\partial_{x_i} k_n) + k_n \partial_{x_i}^2 h\right) W_n - \left((z+h)\partial_{x_i} k_n + k_n \partial_{x_i} h\right)\partial_{x_i} W_n$$
$$+ \left\{\sec^2(k_n H)\left(\partial_{x_i}(k_n H)\right)^2 + \tan(k_n H)\partial_{x_i}^2(k_n H)\right\} Z_n + \tan(k_n H)\partial_{x_i}(k_n H)\partial_{x_i} Z_n.$$
(A20)

To proceed further, we need to calculate the derivative $\partial_{x_i} W_n$. It is calculated similarly as $\partial_{x_i} Z_n$, resulting in

$$\partial_{x_i} W_n = \left((z+h)\partial_{x_i} k_n + k_n \partial_{x_i} h\right) Z_n + \tan(k_n H)\partial_{x_i}(k_n H) W_n. \qquad (A21)$$

Substituting $\partial_{x_i} W_n$ from Eq. (A21), and $\partial_{x_i} Z_n$ from Eq. (A18), into the right-hand side of Eq. (A20), we obtain the second derivative $\partial_{x_i}^2 Z_n$ as a linear combination of $W_n$ and $Z_n$:

$$\partial_{x_i}^2 Z_n = -\left((z+h)\partial_{x_i} k_n + k_n \partial_{x_i} h\right)^2 Z_n$$
$$+ \left\{\sec^2(k_n H)\left(\partial_{x_i}(k_n H)\right)^2 + \tan(k_n H)\partial_{x_i}^2(k_n H) + \left(\tan(k_n H)\cdot\partial_{x_i}(k_n H)\right)^2\right\} Z_n$$
$$- \left((z+h)\partial_{x_i}^2 k_n + 2(\partial_{x_i} h)(\partial_{x_i} k_n) + k_n \partial_{x_i}^2 h\right) W_n$$
$$- 2\left((z+h)\partial_{x_i} k_n + k_n \partial_{x_i} h\right)\tan(k_n H)\partial_{x_i}(k_n H) W_n. \qquad (A22)$$

Now, by using the trigonometric identity $1 + \tan^2 a = \sec^2 a$ and the reference dispersion relation $\tan(k_n H) = -\mu_0/k_n$, Eq. (3.2b), the above equation takes the form

$$\partial_{x_i}^2 Z_n = -\left((z+h)\partial_{x_i} k_n + k_n \partial_{x_i} h\right)^2 Z_n$$
$$+ \left\{\left(\partial_{x_i}(k_n H)\right)^2 - \mu_0 \frac{\partial_{x_i}^2(k_n H)}{k_n} + 2\mu_0^2 \left(\frac{\partial_{x_i}(k_n H)}{k_n}\right)^2\right\} Z_n$$
$$- \left((z+h)\partial_{x_i}^2 k_n + 2(\partial_{x_i} h)(\partial_{x_i} k_n) + k_n \partial_{x_i}^2 h\right) W_n$$
$$+ 2\left((z+h)\partial_{x_i} k_n + k_n \partial_{x_i} h\right)\frac{\mu_0}{k_n}\partial_{x_i}(k_n H) W_n, \qquad (A23)$$

from which we obtain the first equality in Eq. (3.10) after elementary rearrangements of some terms.

- *Detailed proof of the asymptotic estimates appearing in* Eqs. (3.9) – (3.10)

We first state and prove a useful consequence of Proposition 1:

**Lemma A1.** *If $\eta, h \in C^2(\tilde{X}_{\text{part}})$, then the $\mathbf{x}$-derivatives of $k_n H$ satisfy the following asymptotic estimates, uniformly in $\tilde{X}_{\text{part}}$,*

$$\partial_{x_i}(k_n H) = O(n^{-1}), \qquad \partial_{x_i}^2(k_n H) = O(n^{-1}). \tag{A24a,b}$$

*Proof:* Taking the first and the second $x_i$-derivatives of $k_n H$, and substituting $\partial_{x_i} k_n$ from Eq. (3.6) and $\partial_{x_i}^2 k_n$ from Eq. (3.7), we obtain

$$\partial_{x_i}(k_n H) = H \partial_{x_i} k_n + k_n \partial_{x_i} H = \partial_{x_i} H \left( H \partial_H k_n + k_n \right), \tag{A25}$$

$$\partial_{x_i}^2 (k_n H) = H \partial_{x_i}^2 k_n + 2 \partial_{x_i} k_n \partial_{x_i} H + k_n \partial_{x_i}^2 H =$$
$$= (\partial_{x_i} H)^2 \left( H \partial_H^2 k_n + 2 \partial_H k_n \right) + \partial_{x_i}^2 H \left( H \partial_H k_n + k_n \right). \tag{A26}$$

We proceed now to study the terms in parentheses in the right-hand side of Eqs. (A25) and (A26). Using Eq. (A7) (the first equality in Eq. (3.4)), we obtain

$$H \partial_H k_n + k_n = -k_n \left( \frac{H(k_n^2 + \mu_0^2)}{-\mu_0 + H(k_n^2 + \mu_0^2)} - 1 \right) = -k_n \left( \frac{1}{1 - \frac{\mu_0}{H(k_n^2 + \mu_0^2)}} - 1 \right) =$$

$$= -k_n \left( \frac{\mu_0}{H(k_n^2 + \mu_0^2)} + O(k_n^{-4}) \right) = -\frac{\mu_0}{H k_n} + O(k_n^{-3}) = O(n^{-1}), \tag{A27}$$

which proves Eq. (A24a). Further, using Eq. (A10) (the first equality in Eq. (3.5)) and performing simple algebraic manipulations, we obtain

$$H \partial_H^2 k_n + 2 \partial_H k_n = -2 \partial_H k_n \left\{ \mu_0 + \frac{\partial_H k_n}{k_n} (H \mu_0 - 1) \left[ 2 + H \frac{\partial_H k_n}{k_n} \right] \right\} H + 2 \partial_H k_n$$

$$= -2 \partial_H k_n (H \mu_0 - 1) \left( \frac{H \partial_H k_n + k_n}{k_n} \right)^2 = O(n^{-3}). \tag{A28}$$

The last estimate is obtained by using $\partial_H k_n = O(n)$ and Eq. (A13). Combining now Eqs. (A27) and (A28) with Eq. (A26) we find that $\partial_{x_i}^2(k_n H) = O(n^{-1})$, which proves Eq. (A24b). This completes the proof of Lemma A1. ∎

The asymptotic estimates indicated in the second equalities of Eqs. (3.9) - (3.10) are obtained by using the derived analytical expressions for the corresponding quantities, and the estimates

$$k_n = O(n), \quad \partial_{x_i} k_n = O(n), \quad \partial_{x_i}^2 k_n = O(n), \quad \text{(Eqs. (2.5a), (3.6), (3.7))}$$

$$\partial_{x_i}(k_n H) = O(n^{-1}), \quad \partial_{x_i}^2(k_n H) = O(n^{-1}), \quad \text{(Eqs. (A24a,b))}$$

and $Z_n = O(1)$, $W_n = O(1)$. The latter estimates follow easily from the observation

$$\left|Z_n(z)\right|, \left|W_n(z)\right| \leq \frac{1}{\left|\cos(k_n H)\right|} \to 1, \text{ as } n \to \infty, \text{ since } k_n H \to n\pi.$$

Also, we take into account that all geometric, $n-$independent, quantities are of order $O(1)$.

As an example, we give the details of the estimation $\partial_{x_i}^2 Z_n = O(n^2)$. Considering the analytic expression of the second derivative $\partial_{x_i}^2 Z_n$, first equality in Eq. (3.10), and using the above stated estimates, we see that

$$\partial_{x_i}^2 Z_n = -\underbrace{(\partial_{x_i} k_n)^2 (z+h)^2 Z_n}_{O(n^2)} - \underbrace{2k_n(\partial_{x_i} k_n)(\partial_{x_i} h)(z+h) Z_n}_{O(n^2)}$$

$$+ \left\{ -\underbrace{k_n^2(\partial_{x_i} h)^2}_{O(n^2)} + \underbrace{\left(\partial_{x_i}(k_n H)\right)^2}_{O(n^{-2})} - \underbrace{\mu_0 \frac{\partial_{x_i}^2(k_n H)}{k_n}}_{O(n^{-2})} + \underbrace{2\mu_0^2 \left[\frac{\partial_{x_i}(k_n H)}{k_n}\right]^2}_{O(n^{-4})} \right\} \underbrace{Z_n}_{O(1)}$$

$$- \left\{ \underbrace{\partial_{x_i}^2 k_n}_{O(n)} - \underbrace{2\mu_0 \frac{\partial_{x_i} k_n}{k_n} \partial_{x_i}(k_n H)}_{O(n^{-1})} \right\} \underbrace{(z+h) W_n}_{O(1)}$$

$$- \left\{ \underbrace{2(\partial_{x_i} k_n)(\partial_{x_i} h) + k_n \partial_{x_i}^2 h}_{O(n)} - \underbrace{2\mu_0(\partial_{x_i} h)\partial_{x_i}(k_n H)}_{O(n^{-1})} \right\} \underbrace{W_n}_{O(1)} = O(n^2).$$

From the above equation it is clear that the order of $\partial_{x_i}^2 Z_n$ is $O(n^2)$ and this is controlled by first three terms,

$$-(\partial_{x_i} k_n)^2 (z+h)^2 Z_n, \ -2k_n(\partial_{x_i} k_n)(\partial_{x_i} h)(z+h) Z_n \text{ and } -k_n^2(\partial_{x_i} h)^2 Z_n,$$

of the analytic formula for the $\partial_{x_i}^2 Z_n$. All the other terms are of lower order. The proof of Proposition 2 has been now completed.

### Appendix B. Auxiliary results for section 4

Here we assume that the boundary functions $\eta(\boldsymbol{x}), h(\boldsymbol{x})$ are $C^2(\tilde{X}_{\text{part}})$ although some of the results derived below are valid under milder smoothness assumptions.

**(a) Proof of Lemma 1**

- *Proof of Eq.* (4.5).

The estimate $k_n^{-4} = O(n^{-4})$ follows directly from estimate (2.5a), $k_n = O(n)$. For the derivatives of $k_n^{-4}$, we easily find

$$\partial_{x_i} k_n^{-4} = -4 \frac{\partial_{x_i} k_n}{k_n^5} = O(n^{-4}),$$

$$\partial_{x_i}^2 k_n^{-4} = -4 \left[ \frac{\partial_{x_i}^2 k_n}{k_n^5} - 5 \frac{(\partial_{x_i} k_n)^2}{k_n^6} \right] = O(n^{-4}),$$

after exploiting the estimates (3.6) and (3.7) of Proposition 1.

- *Proof of Eq.* (4.6).

The sequence of functions $(\gamma_n(\boldsymbol{x}))_{n \in \mathbb{N}}$ and the corresponding sequences of $\boldsymbol{x}$-derivatives, $(\partial_{x_i} \gamma_n(\boldsymbol{x}))_{n \in \mathbb{N}}$ and $(\partial_{x_i}^2 \gamma_n(\boldsymbol{x}))_{n \in \mathbb{N}}$, are uniformly convergent to finite limits that can be calculated explicitly in terms of $k_n$ and $H = \eta + h$. The existence of these limits (obtained below) trivializes the proof of Eq. (4.6). By direct calculation of the integral involved in the definition of the $L^2$-norm, in conjunction with Eq. (3.2b), we find that $\gamma_n(\boldsymbol{x}) = \|Z_n(\boldsymbol{x}, \cdot)\|_{L^2}^{-2}$ is given by the formula

$$\gamma_n(\boldsymbol{x}) = \sec^{-2}(k_n H) \left[ \frac{H}{2} + \frac{\sin(k_n H) \cos(k_n H)}{2 k_n} \right]^{-1} = \frac{2 k_n^2}{H k_n^2 + f}, \tag{B1}$$

where $f = f(\boldsymbol{x}) = \mu_0^2 H - \mu_0 = O(1)$. Then, by straightforward calculations, we find

$$\partial_{x_i} \gamma_n = \partial_{x_i}(2 k_n^2) \frac{1}{H k_n^2 + f} + 2 k_n^2 \partial_\alpha \left( \frac{1}{H k_n^2 + f} \right) =$$

$$= \frac{4 k_n \partial_{x_i} k_n}{H k_n^2 + f} - \frac{2 k_n^2 \left( k_n^2 \partial_{x_i} H + 2 H k_n \partial_{x_i} k_n + \partial_{x_i} f \right)}{\left( H k_n^2 + f \right)^2}, \tag{B2}$$

$$\partial_{x_i}^2 \gamma_n = \partial_{x_i}^2(2 k_n^2) \frac{1}{H k_n^2 + f} + 2 \partial_{x_i}(2 k_n^2) \partial_{x_i} \left( \frac{1}{H k_n^2 + f} \right) + 2 k_n^2 \partial_{x_i}^2 \left( \frac{1}{H k_n^2 + f} \right) =$$

$$= \frac{4(\partial_{x_i} k_n)^2 + 4 k_n \partial_{x_i}^2 k_n}{H k_n^2 + f} - \frac{8 k_n \partial_{x_i} k_n \left( k_n^2 \partial_{x_i} H + 2 H k_n \partial_{x_i} k_n + \partial_{x_i} f \right)}{\left( H k_n^2 + f \right)^2}$$

$$- \frac{2 k_n^2 \left( k_n^2 \partial_{x_i}^2 H + 4 k_n (\partial_{x_i} H)(\partial_{x_i} k_n) + 2 H (\partial_{x_i} k_n)^2 + 2 H k_n \partial_{x_i}^2 k_n + \partial_{x_i}^2 f \right)}{\left( H k_n^2 + f \right)^2}$$

$$+ \frac{4k_n^2 \left(k_n^2 \partial_{x_i} H + 2H k_n \partial_{x_i} k_n + \partial_{x_i} f\right)^2}{\left(H k_n^2 + f\right)^3}. \tag{B3}$$

Eqs. (B1) - (B3), in conjunction with the fundamental estimate $k_n = \pi n / H + O(n^{-1}) = O(n)$, Eq. (2.5a), and Eqs. (3.6), (3.7) of Proposition 1, permit us to calculate the limits

$$\lim_{n \to \infty} \gamma_n(x) = \frac{2}{H(x)} \stackrel{\text{def}}{=} \gamma(x), \tag{B4}$$

$$\lim_{n \to \infty} \partial_{x_i} \gamma_n(x) = -2 \frac{\partial_{x_i} H(x)}{H^2(x)} = \partial_{x_i} \gamma(x), \tag{B5}$$

$$\lim_{n \to \infty} \partial_{x_i}^2 \gamma_n(x) = -2 \left[ \frac{\partial_{x_i}^2 H(x)}{H^2(x)} - \frac{2 (\partial_{x_i} H(x))^2}{H^3(x)} \right] = \partial_{x_i}^2 \gamma(x). \tag{B6}$$

Note that the convergence of the three sequences is uniform in $\tilde{X}_{\text{part}}$. The fact that the sequences $(\gamma_n(x))_{n \in \mathbb{N}}$, $(\partial_{x_i} \gamma_n(x))_{n \in \mathbb{N}}$ and $(\partial_{x_i}^2 \gamma_n(x))_{n \in \mathbb{N}}$ have well-defined, finite, $n$-independent limits proves Eqs. (4.6). It is worth noting that the commutation between the operations of taking the limit (as $n \to \infty$) and the differentiation, which is apparent in Eqs. (B4) – (B6), is also justified by standard theorems of Calculus. See, for example [1], Theorem 7.17, p. 152.

- *Proof of Eq. (4.7).*

From Eq. (3.3b) we obtain $Z_n(-h) = [Z_n]_{z=-h} = \sec(k_n H)$, which, after differentiations and use of the local dispersion relation (3.2b), leads to

$$\partial_{x_i} Z_n(-h) = -\mu_0 Z_n(-h) \frac{\partial_{x_i}(k_n H)}{k_n}, \tag{B7}$$

$$\partial_{x_i}^2 Z_n(-h) = \mu_0^2 Z_n(-h) \left(\frac{\partial_{x_i}(k_n H)}{k_n}\right)^2 + Z_n^3(-h) \left(\partial_{x_i}(k_n H)\right)^2$$
$$- \mu_0 Z_n(-h) \frac{\partial_{x_i}^2(k_n H)}{k_n}. \tag{B8}$$

The proof of the estimates appearing in Eq. (4.7) is now straightforward. Since $k_n H$ tends to $n\pi$, we get $|\sec(k_n H)| \to 1$, as $n \to \infty$, and thus $Z_n(-h) = O(1)$. Using this result, in conjunction with $k_n = O(n)$, and Eq. (A24a), we find that $\partial_{x_i} Z_n(-h) = O(n^{-2})$. Finally, estimating the terms in Eq. (B8), by using $Z_n(-h) = O(1) = Z_n^3(-h)$, $k_n = O(n)$, and Eqs. (A24a,b), we find $\partial_{x_i}^2 Z_n(-h) = O(n^{-2})$.

**(b) Detailed proof of Lemma 2**

**(i)** Substituting $Z_n = -k_n^{-2}\partial_z^2 Z_n$ (obtained from Eq. (2.3a) with $Q=0$) in the integrand of Eq. (4.8a), performing two integrations by parts in $z$, and using the boundary values of the vertical eigenfunctions ($Z_n(\eta)=1$, $\partial_z Z_n(\eta)=\mu_0$ and $\partial_z Z_n(-h)=0$), we obtain

$$\int_{-h}^{\eta} F Z_n \, dz = -k_n^{-2}\int_{-h}^{\eta} F \partial_z^2 Z_n \, dz =$$

$$= -k_n^{-2}\left\{[F]_{z=\eta}\mu_0 - [\partial_z F]_{z=\eta} + [\partial_z F]_{z=-h}Z_n(-h) + \int_{-h}^{\eta}\partial_z^2 F Z_n \, dz\right\}.$$

Since $F(\boldsymbol{x},\cdot)\in H^2(-h,\eta)$, all boundary values of $F$ and $\partial_z F$ appearing in the right-hand side of the above equation are well-defined and of order $O(1)$. Recall also that $Z_n(-h) = O(1)$. The integral is estimated by using the Cauchy-Schwartz inequality,

$$\left|\int_{-h}^{\eta}\partial_z^2 F Z_n \, dz\right| \le \|\partial_z^2 F\|_{L^2}\|Z_n\|_{L^2} \le \|F\|_{H^2}\|Z_n\|_{L^2} = O(1).$$

Thus, the term in parenthesis is of order $O(1)$, and the order of the integral appearing in the left-hand side of the above equation is $O(k_n^{-2})O(1) = O(n^{-2})$, as asserted.

**(ii)** Substituting $W_n = -k_n^{-2}\partial_z^2 W_n$ in the integrand of Eq. (4.8b), performing two integrations by parts in $z$, and using the boundary values of the auxiliary vertical functions $W_n$,

$$[W_n]_{z=\eta} = \tan(k_n H) = -\frac{\mu_0}{k_n}, \quad [W_n]_{z=-h} = 0,$$

$$[\partial_z W_n]_{z=\eta} = k_n, \quad [\partial_z W_n]_{z=-h} = k_n \sec(k_n H),$$

which are obtained directly form Eq. (3.8), we find

$$\int_{-h}^{\eta} F W_n \, dz = -k_n^{-2}\int_{-h}^{\eta} F \partial_z^2 W_n \, dz =$$

$$= -k_n^{-2}\left\{[F]_{z=\eta}k_n + [\partial_z F]_{z=\eta}\frac{\mu_0}{k_n} - [F]_{z=-h}k_n Z_n(-h) + \int_{-h}^{\eta}\partial_z^2 F W_n \, dz\right\}.$$

Now, the presence of $k_n = O(n)$ in the first and the third terms in parenthesis, makes the latter $O(n)$. Thus, the order of the integral is $O(k_n^{-2})O(n) = O(n^{-1})$, as asserted.

**(iii)** Expressing $\partial_{x_i}Z_n$ in terms of $Z_n$ and $W_n$, from Eq. (3.9), the integral $\int_{-h}^{\eta} F \partial_{x_i} Z_n \, dz$ becomes

$$\int_{-h}^{\eta} F \partial_{x_i} Z_n \, dz =$$
$$= -\partial_{x_i} k_n \int_{-h}^{\eta} F(z+h) W_n \, dz - k_n \partial_{x_i} h \int_{-h}^{\eta} F W_n \, dz - \mu_0 \frac{\partial_{x_i}(k_n H)}{k_n} \int_{-h}^{\eta} F Z_n \, dz.$$

The first integral in the right-hand side of the above equation is $O(n^{-1})$ (apply Eq. (4.8b) with $F_1 = F(z+h)$). Using this result, in conjunction with Eq. (3.6), we see that the first term is $O(1)$. Similarly, the second term is $O(n) \cdot O(n^{-1}) = O(1)$. Finally, using Eq. (A24a), $k_n = O(n)$ and Eq. (4.8a), we deduce that the order of the third term is $O(n^{-1}) \cdot O(n^{-1}) \cdot O(n^{-2}) = O(n^{-4})$. Thus, the order of $\int_{-h}^{\eta} F \partial_{x_i} Z_n \, dz$ is determined by the first term, and it is $O(1)$, which proves the assertion (4.9a).

**(iv)** Now we use Eq. (3.10), expressing $\partial_{x_i}^2 Z_n$ in terms of $Z_n$ and $W_n$. After this substitution, the integral in Eq. (4.9b) takes the form

$$\int_{-h}^{\eta} F \partial_{x_i}^2 Z_n \, dz = -(\partial_{x_i} k_n)^2 \int_{-h}^{\eta} F(z+h)^2 Z_n \, dz - 2 k_n (\partial_{x_i} k_n)(\partial_{x_i} h) \int_{-h}^{\eta} F(z+h) Z_n \, dz$$
$$+ \left\{ -k_n^2 (\partial_{x_i} h)^2 - \left(\partial_{x_i}(k_n H)\right)^2 - \mu_0 \frac{\partial_{x_i}^2 (k_n H)}{k_n} + 2\mu_0^2 \left(\frac{\partial_{x_i}(k_n H)}{k_n}\right)^2 \right\} \int_{-h}^{\eta} F Z_n \, dz$$
$$- \left\{ \partial_{x_i}^2 k_n - 2\mu_0 (\partial_{x_i} k_n) \frac{\partial_{x_i}(k_n H)}{k_n} \right\} \int_{-h}^{\eta} F(z+h) W_n \, dz$$
$$- \left\{ k_n \partial_{x_i}^2 h + 2(\partial_{x_i} k_n)(\partial_{x_i} h) - 2\mu_0 (\partial_{x_i} h) \partial_{x_i}(k_n H) \right\} \int_{-h}^{\eta} F W_n \, dz.$$

All integral terms in the right-hand side of the above equation can be estimated by invoking the already proved assertions (4.8a,b), applied to the functions $F$, $F_1 = F(z+h)$, $F_2 = F(z+h)^2$. Using these results, in conjunction with the known estimates for $k_n$, $\partial_{x_i} k_n$, $\partial_{x_i}^2 k_n$, $\partial_{x_i}(k_n H)$ and $\partial_{x_i}^2 (k_n H)$ (see Eqs. (2.5a), (3.6) and (3.7)) we readily prove (4.9b).

The uniform validity of the estimates (4.8) and (4.9), asserted in the last part of the proposition, follows from the fact that the uniform validity of the assumptions ensures that all bounds (order relations) apply uniformly for all $x \in \tilde{X}_{part}$.

**Remark.** The four integrals in Eqs. (4.8) and (4.9) could be directly estimated under the weaker assumption $F \in L^2(-h, \eta)$, using the Cauchy-Schwartz inequality. However, the estimates obtained in this way are weaker, namely **(i):** $O(1)$, **(ii):** $O(1)$, **(iii):** $O(n)$, and **(iv):** $O(n^2)$, which are not satisfactory for our purposes. Assuming that $F$ is sufficiently smooth, i.e. $F \in H^2(-h, \eta)$, we are able to perform two integrations by parts, leading to the better estimates derived above.

### (c) Detailed proof of Lemma 3

The proof of $\lambda_n(x) = O(1)$ is fully described in the main part of the paper. To prove $\partial_{x_i} \lambda_n = O(1)$, we differentiate Eq. (4.4) with respect to $x_i$, obtaining

$$\partial_{x_i} \lambda_n = \partial_{x_i}\left(\left[\partial_z^2 \Phi^* \mu_0 - \partial_z^3 \Phi^*\right]_{z=\eta}\right) - \partial_{x_i}\left(\left[\partial_z^3 \Phi^*\right]_{z=-h} Z_n(-h)\right)$$
$$+ \partial_{x_i} \eta \left[\partial_z^4 \Phi^*\right]_{z=\eta} + \partial_{x_i} h \left[\partial_z^4 \Phi^*\right]_{z=-h} Z_n(-h) +$$
$$+ \int_{-h}^{\eta} \partial_{x_i} \partial_z^4 \Phi^* Z_n \, dz + \int_{-h}^{\eta} \partial_z^4 \Phi^* \partial_{x_i} Z_n \, dz. \quad (B9)$$

Smoothness assumptions and the estimates (4.7) ensure that all boundary terms (the first two lines of the right-hand side of Eq. (B9)) are of order $O(1)$. Applying Eqs. (4.8a) and (4.9a), with $F = \partial_{x_i} \partial_z^4 \Phi^* \in H^2(-h, \eta)$ and $F = \partial_z^4 \Phi^* \in H^2(-h, \eta)$, respectively, it is easily seen that the integral terms in the third line of (B9) are also $O(1)$. Having established that all terms in the right-hand side of Eq. (B9) are of order $O(1)$, we conclude that $\partial_{x_i} \lambda_n = O(1)$.

To prove $\partial_{x_i}^2 \lambda_n = O(1)$, we need to calculate $\partial_{x_i}^2 \lambda_n$. Differentiating Eq. (B9) with respect to $x_i$, we obtain

$$\partial_{x_i}^2 \lambda_n = \partial_{x_i}^2 \left[\partial_z^2 \Phi^* \mu_0 - \partial_z^3 \Phi^*\right]_{z=\eta} - \partial_{x_i}^2 \left(\left[\partial_z^3 \Phi^*\right]_{z=-h} Z_n(-h)\right)$$
$$+ \partial_{x_i}\left(\partial_{x_i} \eta \left[\partial_z^4 \Phi^*\right]_{z=\eta}\right) + \partial_{x_i}\left(\partial_{x_i} h \left[\partial_z^4 \Phi^*\right]_{z=-h} Z_n(-h)\right)$$
$$+ \partial_{x_i} \eta \left[\partial_{x_i} \partial_z^4 \Phi^* + \partial_z^4 \Phi^* \partial_{x_i} Z_n\right]_{z=\eta} + \partial_{x_i} h \left[\partial_{x_i} \partial_z^4 \Phi^* + \partial_z^4 \Phi^* \partial_x Z_n\right]_{z=-h} Z_n(-h)$$
$$+ \int_{-h}^{\eta} \partial_{x_i}^2 \partial_z^4 \Phi^* Z_n \, dz + 2\int_{-h}^{\eta} \partial_{x_i} \partial_z^4 \Phi^* \partial_{x_i} Z_n \, dz + \int_{-h}^{\eta} \partial_z^4 \Phi^* \partial_{x_i}^2 Z_n \, dz. \quad (B10)$$

Again, under the stated smoothness assumptions and Eqs. (4.7), all boundary terms (the first three lines of the right-hand side of Eq. (B10)) are of order $O(1)$. The integrals, appearing in the last line of the right-hand side of Eq. (B10), are estimated by applying Eqs. (4.8a), (4.9a) and (4.9b) with $F = \partial_{x_i}^2 \partial_z^4 \Phi^*$, $F = \partial_{x_i} \partial_z^4 \Phi^*$ and $F = \partial_z^4 \Phi^*$, respectively, and are found to be of order $O(1)$. These findings ensures that $\partial_{x_i}^2 \lambda_n = O(1)$. Thus, the proof of Lemma 2 has been completed.

**Appendix C. Variational derivation of the CCMS, Eqs.** (5.2), (5.3), **and the Dirichlet-to-Neumann operator, Eq.** (5.5)

**(a) Variational derivation of Eqs.** (5.2), (5.3)

In this appendix we exploit the exact semi-separation of variables, $\Phi = \sum_{n=-2}^{\infty} \varphi_n Z_n$, which, in accordance with the results of section 4, can be termwise differentiated two times throughout the closed domain $\bar{D}_h^\eta(X)$, in order to derive exact modal reformulations of the following boundary-value problem (BVP):

$$\Delta \Phi = 0, \quad \text{in} \quad D_h^\eta(X), \tag{C1a}$$

$$\left[\Phi\right]_{z=\eta} = \psi, \qquad \left[\partial_n \Phi\right]_{z=-h} = 0, \tag{C1b,c}$$

$$\left[\partial_n \Phi\right]_{x_1=a} = V_a, \quad \nabla\Phi \to 0 \text{ as } |\boldsymbol{x}| \to \infty. \tag{C1c,d}$$

Although the laterally periodic case is examined in the main paper, here we choose to work with a laterally semi-bounded domain with boundary conditions (C1c,d) in order to show the versatility of the present approach, concerning the lateral boundary conditions. The derivation that follows is easily adapted to other cases, and the statement of the corresponding results can be found in subsection (c) below.

Our starting point for the derivation of Eqs. (5.2) and (5.3) is the following variational formulation of the BVP (C1) (see, e.g., [2]):

$$\delta \mathcal{I}[\Phi, \lambda; \delta\Phi, \delta\lambda] \equiv \delta_\Phi \mathcal{I}[\Phi, \lambda; \delta\Phi] + \delta_\lambda \mathcal{I}[\Phi, \lambda; \delta\lambda] = 0, \tag{C2a}$$

where the functional $\mathcal{I}[\Phi, \lambda]$ is given by

$$\mathcal{I}[\Phi, \lambda] = \frac{1}{2} \int_X \int_{-h}^{\eta} |\nabla\Phi|^2 \, dz \, d\boldsymbol{x} - \int_X \lambda \left(\left[\Phi\right]_{z=\eta} - \psi\right) d\boldsymbol{x} \\ + \int_{\mathbb{R}} \int_{-h_a}^{\eta_a} V_a \left[\Phi\right]_{x_1=a} dz \, dx_2. \tag{C2b}$$

Note that, in the above variational formulation, Dirichlet boundary condition (C1b) is introduced into the functional by means of the Lagrange multiplier $\lambda = \lambda(\boldsymbol{x}, t)$, leaving $\Phi$ unconstraint throughout the closed domain $\bar{D}_h^\eta(X)$, apart from the requirements (essential condition) $\nabla\Phi \to 0$ as $|\boldsymbol{x}| \to \infty$. Calculating the first variations $\delta_\Phi \mathcal{I}[\Phi, \lambda; \delta\Phi]$ and $\delta_\lambda \mathcal{I}[\Phi, \lambda; \delta\lambda]$ of the functional $\mathcal{I}[\Phi, \lambda]$ in the usual way, we find

$$\delta_\Phi \mathcal{I}[\Phi, \lambda; \delta\Phi] = \int_X \left\{ \int_{-h}^{\eta} (\nabla\Phi)\cdot(\nabla\delta\Phi) \, dz - \lambda [\delta\Phi]_{z=\eta} \right\} d\boldsymbol{x} + \int_{\mathbb{R}} \int_{-h_a}^{\eta_a} V_a [\delta\Phi]_{x_1=a} \, dz \, dx_2$$

$$= -\int_X \int_{-h}^{\eta} \Delta\Phi \, \delta\Phi \, dz \, d\boldsymbol{x} - \int_X \left(\lambda - \boldsymbol{N}_\eta \cdot [\nabla\Phi]_{z=\eta}\right) [\delta\Phi]_{z=\eta} \, d\boldsymbol{x} +$$

$$+ \int_X \left(\boldsymbol{N}_h \cdot [\nabla\Phi]_{z=-h}\right) [\delta\Phi]_{z=-h} \, d\boldsymbol{x} + \int_{\mathbb{R}} \int_{-h_a}^{\eta_a} \left(-[\partial_{x_1}\Phi]_{x_1=a} + V_a\right) [\delta\Phi]_{x_1=a} \, dz \, dx_2, \tag{C3a}$$

$$\delta_\lambda \mathcal{I}[\Phi, \lambda; \delta\lambda] = -\int_X \left([\Phi]_{z=\eta} - \psi\right) \delta\lambda \, d\boldsymbol{x}. \tag{C3b}$$

In accordance to Eqs. (C3a,b), the two partial variational equations $\delta \mathcal{I}_\Phi[\Phi, \lambda; \delta\Phi] = 0$ and $\delta_\lambda \mathcal{I}[\Phi, \lambda; \delta\lambda] = 0$, in conjunction with the arbitrariness of $\delta\Phi$ and $\delta\lambda$, proves that the variational principle (C2) is equivalent to the BVP (C1) and the equation

$$\lambda = \boldsymbol{N}_\eta \cdot [\nabla \Phi]_{z=\eta}, \tag{C4}$$

where $\boldsymbol{N}_\eta = (-\nabla_{\boldsymbol{x}} \eta, 1)$ is the outward (with respect to the domain $D_h^\eta$) normal (but not normalized to have unit length) vector on the upper boundary $\Gamma^\eta$.

We shall now apply the variational formulation (C2) in conjunction with the series representation of the wave potential $\Phi$, in order to prove Eqs. (5.2) and (5.3). By the substitution $\Phi(\boldsymbol{x}, z, t) = \sum_{n=-2}^{\infty} \varphi_n Z_n \equiv \boldsymbol{\varphi}^T \boldsymbol{Z}$, where $\boldsymbol{Z} = \boldsymbol{Z}(z, \eta, h) \equiv \{Z_n(z; \eta(\boldsymbol{x},t), h(\boldsymbol{x}))\}_{n \geq -2}$ is the system of vertical basis functions given by Eqs. (2.8a,b), (3.3a,b), and $\boldsymbol{\varphi} \equiv \boldsymbol{\varphi}(\boldsymbol{x}, t) = \{\varphi_n(\boldsymbol{x}, t)\}_{n \geq -2}$ is the corresponding system of the unknown modal amplitudes, we introduce the change of (functional) variables $(\Phi, \lambda) \to (\boldsymbol{\varphi}, \lambda)$ in the functional (C2b), obtaining the transformed functional $\tilde{\mathcal{I}}[\boldsymbol{\varphi}, \lambda] \equiv \mathcal{I}[\boldsymbol{\varphi}^T \boldsymbol{Z}, \lambda]$. Note that, under the assumptions stated in Theorem 2, termwise differentiation and integration of the infinite series $\boldsymbol{\varphi}^T \boldsymbol{Z}$ is permissible (Corollary 2). Thus, it is legitimate to perform termwise manipulations in the right-hand side of Eq. (C2b), after the substitution $\Phi(\boldsymbol{x}, z, t) = \boldsymbol{\varphi}^T \boldsymbol{Z}$. Further, under the same assumptions, the functional transformation $(\Phi, \lambda) \to (\boldsymbol{\varphi}, \lambda)$ is continuous, differentiable and invertible, which ensures that the variational equation (C2a) is equivalent[14] to the following, transformed variational equation:

$$\delta \tilde{\mathcal{I}}[\boldsymbol{\varphi}, \lambda; \delta\boldsymbol{\varphi}, \delta\lambda] \equiv \delta_\lambda \tilde{\mathcal{I}}[\boldsymbol{\varphi}, \lambda; \delta\lambda] + \sum_{m=-2}^{\infty} \delta_{\varphi_m} \tilde{\mathcal{I}}[\boldsymbol{\varphi}, \lambda; \delta\varphi_m] = 0. \tag{C5}$$

Performing the variations indicated in Eq. (C5), we easily obtain the following Euler-Lagrange equations

$$\delta\varphi_m: \quad \int_{-h}^{\eta} \Delta(\boldsymbol{\varphi}^T \boldsymbol{Z}) Z_m \, dz - \boldsymbol{N}_h \cdot [\nabla(\boldsymbol{\varphi}^T \boldsymbol{Z}) Z_m]_{z=-h} = \\ = -\lambda + \boldsymbol{N}_\eta \cdot [\nabla(\boldsymbol{\varphi}^T \boldsymbol{Z})]_{z=\eta}, \quad m \geq -2, \; \boldsymbol{x} \in X, \tag{C6a}$$

$$[\delta\varphi_m]_{x_1=a}: \quad \int_{-h_a}^{\eta_a} \left(-[\partial_{x_1}(\boldsymbol{\varphi}^T \boldsymbol{Z})]_{x_1=a} + V_a\right) [Z_m]_{x_1=a} \, dz = 0, \; m \geq -2, \; x_2 \in \mathbb{R}, \tag{C6b}$$

$$\delta\lambda: \quad [\boldsymbol{\varphi}^T \boldsymbol{Z}]_{z=\eta} = \psi. \tag{C6c}$$

---

[14] That is, $(\Phi, \lambda)$ solves Eq. (C2a) if and only if $(\boldsymbol{\varphi}, \lambda)$ solves Eq. (C5).

Observe now that the right-hand side of Eq. (C6a) is zero because of Eq. (C4). Accordingly, Eq. (C6a) can be written as

$$\int_{-h}^{\eta} \Delta(\boldsymbol{\varphi}^T \boldsymbol{Z}) Z_m \, dz - \boldsymbol{N}_h \cdot \left[ \nabla(\boldsymbol{\varphi}^T \boldsymbol{Z}) Z_m \right]_{z=-h} = 0, \quad m \geq -2, \quad \boldsymbol{x} \in X, \tag{C7}$$

which, after substituting $\boldsymbol{\varphi}^T \boldsymbol{Z} = \sum_{n=-2}^{\infty} \varphi_n Z_n$, performing (termwise) the indicated differentiation and integration, and using the notation $\Delta = \nabla_x^2 + \partial_z^2$, $\nabla = (\nabla_x, \partial_z)$ and $\boldsymbol{N}_h = (-\nabla_x h, -1)$, takes the form

$$\sum_{n=-2}^{\infty} \int_{-h}^{\eta} \left( (\nabla_x^2 \varphi_n) Z_n + 2(\nabla_x \varphi_n) \cdot (\nabla_x Z_n) + \varphi_n (\nabla_x^2 Z_n + \partial_z^2 Z_n) \right) Z_m \, dz \; +$$

$$+ \sum_{n=-2}^{\infty} \left\{ (\nabla_x h) \cdot \left( (\nabla_x \varphi_n) \left[ Z_n Z_m \right]_{z=-h} + \varphi_n \left[ (\nabla_x Z_n) Z_m \right]_{z=-h} \right) + \varphi_n \left[ \partial_z Z_n Z_m \right]_{z=-h} \right\} = 0. \tag{C8a}$$

Note that the last two terms in the second series can be recombined in the form

$$- \varphi_n \boldsymbol{N}_h \cdot \left[ (\nabla_x Z_n, \partial_z Z_n) Z_m \right]_{z=-h},$$

and thus the second series in the left-hand side of Eq. (C8a) takes the form

$$\sum_{n=-2}^{\infty} \left\{ (\nabla_x \varphi_n) \cdot (\nabla_x h) \left[ Z_n Z_m \right]_{z=-h} - \varphi_n \boldsymbol{N}_h \cdot \left[ (\nabla_x Z_n, \partial_z Z_n) Z_m \right]_{z=-h} \right\}. \tag{C8b}$$

Observing now that $\varphi_n$ and their derivatives are independent from $z$ (thus, they can be taken out of the vertical integrals), taking into account (C8b), and collecting together the terms containing $\nabla_x^2 \varphi_n$, $\nabla_x \varphi_n$ and $\varphi_n$, we conclude that Eq. (C8a) leads to the Eq. (5.2a) and Eqs. (5.3a-c). Working similarly, we find that Eq. (C6b), results in

$$\sum_{n=-2}^{\infty} T_{mn}^{(a)} \varphi_n (a, x_2) = \int_{-h_a}^{\eta_a} V_a \left[ Z_m \right]_{x_1=a} dz \equiv g_m(x_2), \quad m \geq -2, \; x_2 \in \mathbb{R}, \tag{C9}$$

where the boundary operators $T_{mn}^{(a)}$ are given by

$$T_{mn}^{(a)} = \left[ A_{mn} \partial_{x_1} + (1/2) \left( B_{mn}^{(1)} - \partial_{x_1} h \left[ Z_m Z_n \right]_{z=-h} \right) \right]_{x_1 = a}.$$

Eq. (C9) is the modal form of the lateral boundary condition (C1c). Finally, recalling that $\left[ Z_n \right]_{z=\eta} = 1$, $n \geq -2$, we see that Eq.(C6c) results in

$$\sum_{n=-2}^{\infty} \varphi_n = \psi,$$

which is identical with Eq. (5.2b).

### (b) The form of the corresponding DtN operator, in terms of the enhanced modal series expansion

It remains to prove the modal form of the DtN operator, Eq. (5.5). Observe that Eq. (C4) provides a physical interpretation of the Lagrange multiplier $\lambda = \lambda(\boldsymbol{x},t)$; it is a representation of the DtN operator (5.4), that is, $\lambda(\boldsymbol{x},t) = \mathsf{G}[\eta, h, V_a]\psi(\boldsymbol{x},t)$. Substituting the modal expansion $\Phi = \sum_{n=-2}^{\infty} \varphi_n Z_n$ into the right-hand side of Eq. (C4), and taking into account the condition $[Z_n]_{z=\eta} = 1$, we find

$$\mathsf{G}[\eta, h, V_a]\psi = -\nabla_x \eta \cdot \sum_{n=-2}^{\infty} \left(\nabla_x \varphi_n + \varphi_n [\nabla_x Z_n]_{z=\eta}\right) + \sum_{n=-2}^{\infty} \varphi_n [\partial_z Z_n]_{z=\eta}. \quad \text{(C10)}$$

Eq. (C10) is greatly simplified by using the analytic expressions for the boundary values of the derivatives $[\nabla_x Z_n]_{z=\eta}$ and $[\partial_z Z_n]_{z=\eta}$, $n \geq -2$. Such expressions have been already found in section 3 for $n \geq 1$; see Eq. (3.9). A similar result holds true for $n = 0$ as well, leading to

$$[\nabla Z_n]_{z=\eta} = \begin{pmatrix} -\mu_0 \nabla_x \eta \\ \mu_0 \end{pmatrix}, \quad n \geq 0. \quad \text{(C11a)}$$

For $n = -2, -1$, the corresponding formulae are obtained by direct differentiation of Eqs. (2.8a,b):

$$[\nabla Z_{-2}]_{z=\eta} = \begin{pmatrix} -(\mu_0 + h_0^{-1})\nabla_x \eta \\ \mu_0 + h_0^{-1} \end{pmatrix}, \quad [\nabla Z_{-1}]_{z=\eta} = \begin{pmatrix} -\mu_0 \nabla_x \eta \\ \mu_0 \end{pmatrix}. \quad \text{(C11b)}$$

Substituting Eq. (C11a,b) into Eq. (C10), and interchanging differentiations with the infinite summation, which is permissible under the considered smoothness assumptions (Corollary 2, section 4), we obtain

$$\mathsf{G}[\eta, h, V_a]\psi = -\nabla_x \eta \cdot \nabla_x \sum_{n=-2}^{\infty} \varphi_n +$$
$$+ (\nabla_x \eta)^2 \left(\mu_0 \sum_{n=-2}^{\infty} \varphi_n + h_0^{-1} \varphi_{-2}\right) + \mu_0 \sum_{n=-2}^{\infty} \varphi_n + h_0^{-1} \varphi_{-2}. \quad \text{(C12)}$$

Finally, using Eq. (C9), Eq. (C12) becomes

$$\mathsf{G}[\eta, h]\psi = -(\nabla_x \eta) \cdot (\nabla_x \psi) + (|\nabla_x \eta|^2 + 1)\left(h_0^{-1} \varphi_{-2} + \mu_0 \psi\right),$$

which is the same as Eq. (5.5).

### (c) Complete formulation of some specific boundary-value problems as CCMS

We now present the forms of the CCMS for four geometrical configurations different than the one considered above.

**The 2D periodic problem;** $X = \mathbb{R}$

In this case, we assume that the functions $\eta = \eta(x,t), h = h(x)$ and $\psi = \psi(x)$ are periodic functions of $x \in X$, with period (say) $2\pi$. Then, the modal form of the BVP (5.1) is as follows:

$$\sum_{n=-2}^{\infty} L_{mn}[\eta,h]\varphi_n \equiv \sum_{n=-2}^{\infty} \left(A_{mn}\partial_x^2 + B_{mn}\partial_x + C_{mn}\right)\varphi_n = 0, \quad m \geq -2, \quad x \in [0, 2\pi], \quad \text{(C13a)}$$

$$\sum_{n=-2}^{\infty} \varphi_n = \psi, \quad x \in [0, 2\pi], \quad \text{(C13b)}$$

$$\varphi_n(x) = \varphi_n(x + 2\pi), \quad \text{(C13c)}$$

where

$$A_{mn} = \int_{-h}^{\eta} Z_n Z_m \, dz, \quad \text{(C13d)}$$

$$B_{mn} = 2\int_{-h}^{\eta} (\partial_x Z_n) Z_m \, dz + (\partial_x h)\left[Z_m Z_n\right]_{z=-h}, \quad \text{(C13e)}$$

$$C_{mn} = \int_{-h}^{\eta} (\partial_x^2 Z_n + \partial_z^2 Z_n) Z_m \, dz - N_h \cdot \left[(\partial_x Z_n, \partial_z Z_n) Z_m\right]_{z=-h}. \quad \text{(C13f)}$$

and $\quad N_h = (-\partial_x h, -1).$ (C13g)

**The 3D periodic problem;** $X = \mathbb{R}^2$

In this case, we assume that the functions $\eta = \eta(\boldsymbol{x},t), h = h(\boldsymbol{x})$ and $\psi = \psi(\boldsymbol{x})$, where $\boldsymbol{x} = (x_1, x_2)$, are doubly-periodic functions of $\boldsymbol{x} \in X$, with period (say) $(2\pi, 2\pi)$. Then, the CCMS takes the following form:

$$\text{Eq. (5.2a), with } \boldsymbol{x} \in [0, 2\pi] \times [0, 2\pi], \quad \text{(C14a)}$$

$$\sum_{n=-2}^{\infty} \varphi_n = \psi, \quad \boldsymbol{x} \in [0, 2\pi] \times [0, 2\pi], \quad \text{(C14b)}$$

$$\varphi_n(x_1, x_2) = \varphi_n(x_1 + 2\pi, x_2 + 2\pi), \quad \text{(C14c)}$$

with $A_{mn}, B_{mn}, C_{mn}$ given by Eqs. (5.3a,b,c).

**2D closed basin with vertical lateral boundaries;** $X = [a,b]$

$$\sum_{n=-2}^{\infty} L_{mn}[\eta,h]\varphi_n \equiv \sum_{n=-2}^{\infty} \left(A_{mn}\partial_x^2 + B_{mn}\partial_x + C_{mn}\right)\varphi_n = 0, \quad m \geq -2, \quad x \in [a,b], \quad \text{(C15a)}$$

$$\sum_{n=-2}^{\infty} \varphi_n = \psi, \quad x \in [a,b], \quad \text{(C15b)}$$

$$\left[\sum_{n=-2}^{M}\left(A_{mn}\partial_x\varphi_n + \frac{1}{2}\Big(B_{mn} - \partial_x h\left[Z_m Z_n\right]_{z=-h}\Big)\varphi_n\right)\right]\Bigg|_{x=a} = 0, \tag{C15c}$$

$$\left[\sum_{n=-2}^{M}\left(A_{mn}\partial_x\varphi_n + \frac{1}{2}\Big(B_{mn}\varphi_n - \partial_x h\left[Z_m Z_n\right]_{z=-h}\Big)\right)\right]\Bigg|_{x=b} = 0, \tag{C15d}$$

with $A_{mn}$, $B_{mn}$, $C_{mn}$ given by Eqs. (C13d,e,f) and $N_h$ as in Eq. (C13g).

**3D Closed basin with vertical lateral boundaries; $X$ is a bounded domain in $\mathbb{R}^2$, with boundary $\partial X$ defined by the closed curve $\Gamma$: $f(x_1, x_2) = 0$**

$$\text{Eqs. (5.2a), with } \boldsymbol{x} \in X, \tag{C16a}$$

$$\sum_{n=-2}^{\infty}\varphi_n = \psi, \quad \boldsymbol{x} \in X, \tag{C16b}$$

$$\sum_{n=-2}^{\infty}\mathrm{T}_{mn}^{(\partial X)}\varphi_n(\boldsymbol{x}) = 0, \quad m \geq -2, \quad \boldsymbol{x} \in \partial X, \tag{C16c}$$

with $A_{mn}$, $\boldsymbol{B}_{mn}$, $C_{mn}$ given by Eqs. (5.3a,b,c),

$$\mathrm{T}_{mn}^{(\partial X)} = \boldsymbol{n}_\Gamma \cdot \left[A_{mn}\nabla_{\boldsymbol{x}} + \frac{1}{2}\Big(\boldsymbol{B}_{mn} - \nabla_{\boldsymbol{x}} h\left[Z_m Z_n\right]_{z=-h}\Big)\right]\Bigg|_{\boldsymbol{x}\in\partial X}, \tag{C16d}$$

and $\boldsymbol{n}_\Gamma = -\nabla_{\boldsymbol{x}} f(\boldsymbol{x}) / \|\nabla_{\boldsymbol{x}} f(\boldsymbol{x})\|$ is the outward unit normal vector on the lateral vertical wall defined by $\Gamma$.

## Appendix D. Closed-form expressions for the modal amplitudes $\{\varphi_n(x)\}_{n=-2}^{\infty}$ of the harmonic function $\Phi_\kappa(x, z)$

According to Theorem 1, any sufficiently smooth function $\Phi(x, z)$, defined on $\bar{D}_h^\eta(X)$, can be expanded in the form $\Phi(x, z) = \sum_{n=-2}^{\infty}\varphi_n(x)\, Z_n(z; \eta(x), h(x))$, where the modal amplitudes $\{\varphi_n(x)\}_{n=-2}^{\infty}$ are defined by Eqs. (2.9) in terms of $\Phi$. Consider now the $x$-periodic function $\Phi(x, z) = \Phi_\kappa(x, z) = \cosh(\kappa(z + h_0))\cos(\kappa x)$ that satisfies Eqs. (5.1) with $\psi = \cosh(\kappa(\eta + h_0))\cos(\kappa x)$. Its modal amplitudes $\{\varphi_n(x)\}_{n=-2}^{\infty}$ can be explicitly calculated in view of Eqs. (2.9). Indeed, Eqs. (2.9a,b) give

$$\varphi_{-2} = h_0\Big(\kappa \sinh(\kappa H) - \mu_0 \cosh(\kappa H)\Big)\cos(\kappa x), \tag{D1}$$

$$\varphi_{-1} = 0, \tag{D2}$$

where $H = \eta + h$. Then, from Eq. (2.9c) it follows that

$$\Phi^* = \cosh(\kappa(z + h_0))\cos(\kappa x) - \varphi_{-2} Z_{-2},$$

and from Eq. (2.9d) we obtain (for $n \geq 0$)

$$\varphi_n = \|Z_n\|_{L^2}^{-2}\left[\left(\int_{-h}^{\eta}\cosh(\kappa(z+h_0))Z_n\,dz\right)\cos(\kappa x) - \varphi_{-2}\int_{-h}^{\eta}Z_{-2}Z_n\,dz\right]. \quad (D3)$$

For $n = 0$, the first integral in the right-hand side of the above equation becomes

$$\int_{-h}^{\eta}\cosh(\kappa(z+h_0))\frac{\cosh(k_0(z+h_0))}{\cosh(k_0 H)}\,dz =$$

$$= \frac{1}{2\cosh(k_0 H)}\int_{-h}^{\eta}\{\cosh[(\kappa+k_0)(z+h_0)] + \cosh[(\kappa-k_0)(z+h_0)]\}dz$$

$$= \frac{1}{2\cosh(k_0 H)}\left(\frac{\sinh[(\kappa+k_0)H]}{\kappa+k_0} + \frac{\sinh[(\kappa-k_0)H]}{\kappa-k_0}\right)$$

$$= \frac{\kappa\sinh(\kappa H) - k_0\cosh(\kappa H)\tanh(k_0 H)}{\kappa^2 - k_0^2}$$

$$= \frac{\kappa\sinh(\kappa H) - \mu_0\cosh(\kappa H)}{\kappa^2 - k_0^2},$$

where Eq. (3.2a), $\tanh(k_0 H) = \mu_0 / k_0$, has been used in the last equality. Using this result in conjunction Eq. (D1), we obtain that

$$\left(\int_{-h}^{\eta}\cosh(\kappa(z+h_0))Z_n\,dz\right)\cos(\kappa x) = \frac{\varphi_{-2}}{h_0(\kappa^2 - k_0^2)}. \quad (D4)$$

The second integral in the right-hand side of Eq. (D3), for $n = 0$, is calculated as follows:

$$\int_{-h}^{\eta}Z_{-2}Z_0\,dz = \int_{-h}^{\eta}\left(\alpha\frac{(z+h)^2}{H} - \alpha H + 1\right)\frac{\cosh(k_0(z+h_0))}{\cosh(k_0 H)}\,dz$$

$$= \alpha\left(-\frac{2}{k_0^2} + \frac{2\tanh(k_0 H)}{k_0^3 H} + \frac{H\tanh(k_0 H)}{k_0}\right) + (-\alpha H + 1)\frac{\tanh(k_0 H)}{k_0}$$

$$= 2\alpha\left(-\frac{1}{k_0^2} + \frac{\mu_0}{k_0^4 H}\right) + \frac{\mu_0}{k_0^2} = 2\alpha\mu_0\frac{1}{k_0^4 H} + (\mu_0 - 2a)\frac{1}{k_0^2}, \quad (D5)$$

where $\alpha = (\mu_0 h_0 + 1)/2h_0$ and, once again, Eq. (3.2a) has been used in the third equality. Substituting Eqs. (D4) and (D5) into Eq. (D3), for $n = 0$, we find

$$\varphi_0 = \frac{\varphi_{-2}}{\|Z_0\|_{L^2}^2}\left(\frac{1}{h_0(\kappa^2 - k_0^2)} - \frac{\mu_0(\mu_0 h_0 + 1)}{h_0 H k_0^4} - \left(\mu_0 - \frac{\mu_0 h_0 + 1}{h_0}\right)\frac{1}{k_0^2}\right) =$$

$$= \frac{\varphi_{-2}}{\|Z_0\|_{L^2}^2}\left(\frac{k_0^2 - \mu_0 h_0(\kappa^2 - k_0^2) + (\mu_0 h_0 + 1)(\kappa^2 - k_0^2)}{h_0(\kappa^2 - k_0^2)k_0^2} - \frac{\mu_0(\mu_0 h_0 + 1)}{h_0 H k_0^4}\right) =$$

$$= \frac{\varphi_{-2}}{\|Z_0\|_{L^2}^2} \left( \frac{\kappa^2}{h_0(\kappa^2 - k_0^2)k_0^2} - \frac{\mu_0(\mu_0 h_0 + 1)}{h_0 H k_0^4} \right). \tag{D6}$$

Similarly, one can calculate

$$\varphi_n = \frac{\varphi_{-2}}{\|Z_n\|_{L^2}^2} \left( -\frac{\kappa^2}{h_0(\kappa^2 + k_n^2)k_n^2} - \frac{\mu_0(\mu_0 h_0 + 1)}{h_0 H k_n^4} \right), \qquad n \geq 1. \tag{D7}$$

In the above equations, the norms $\|Z_0\|_{L^2}^2$ and $\|Z_n\|_{L^2}^2$ are given by the equations

$$\|Z_0\|_{L^2}^2 = \frac{H(k_0^2 - \mu_0^2) + \mu_0}{2k_0^2}, \qquad \|Z_n\|_{L^2}^2 = \frac{H(k_n^2 + \mu_0^2) - \mu_0}{2k_n^2}, \quad n \geq 1. \tag{D8}$$

Eqs. (D1), (D2), (D6) and (D7) provide us with closed-form expressions for the modal amplitudes $\{\varphi_n(x)\}_{n=-2}^{\infty}$ of the function $\Phi_\kappa(x,z)$, defined on the non-uniform strip-like domain $0 \leq x \leq 2\pi$, $-h_0 \leq z \leq \eta(x)$, for any surface profile $\eta(x)$. Note that all $\varphi_n$ are expressed in terms of the eigenvalues $k_n = k_n(x)$, $n \geq 0$, and the known quantities $\mu_0, h_0, \kappa$ and $H = H(x,t)$. Two remarks are in order here:

**Remarks: D1.** The eigenvalues $k_n = k_n(x)$ are calculated by solving numerically the transcendental equations

$$\mu_0 - k_0(x) \tanh\big(k_0(x)(\eta(x) + h_0)\big) = 0,$$
$$\mu_0 + k_n(x) \tan\big(k_n(x)(\eta(x) + h_0)\big) = 0, \qquad n \geq 1,$$

at machine accuracy, e.g., by the Newton-Raphson method.

**D2.** It can be shown that the asymptotic behaviour of $\varphi_n(x)$, given by Eqs. (D7), is $\varphi_n(x) = O(n^{-4})$, in accordance with the general estimate obtained in Theorem 1. Moreover, by calculating the derivatives $\partial_{x_i} \varphi_n(x)$ and $\partial_{x_i}^2 \varphi_n(x)$ and using the results of Proposition 1 concerning $\partial_{x_i} k_n(x)$ and $\partial_{x_i}^2 k_n(x)$ we can also verify the decay rates $\partial_{x_i} \varphi_n(x) = O(n^{-4})$ and $\partial_{x_i}^2 \varphi_n(x) = O(n^{-4})$, in accordance with Theorem 2.

## Appendix E. Supplementary numerical results and discussion

### (a) Computation of the DtN operator. Convergence and efficiency of the modal approach

The convergence and the accuracy of the numerical approximation of the DtN operator, $G^{(N_{tot})}[\eta, h_0]\psi$, can be investigated by comparison with its exact values $\mathcal{G}[\eta, h_0]\psi$, explicitly calculated by plugging $\Phi_\kappa(x,z)$ and $\eta(x)$ into Eq. (5.4) [3]. Results for the $L^2$-relative error

$$ER^{(N_{tot})}[\mathcal{G}] := \left\| \mathcal{G}[\eta, h_0]\psi - \mathcal{G}^{(N_{tot})}[\eta, h_0]\psi \right\|_{L^2(X)} / \left\| \mathcal{G}[\eta, h_0]\psi \right\|_{L^2(X)}$$

have been obtained for $\varepsilon = 0.1(0.2)0.9$, by using 3 up to 70 modes, and are shown in **Figure E1**. It is readily observed that: **1)** In both the smooth and the rough case $ER^{(N_{tot})}[\mathcal{G}]$ decays rapidly, up to a plateau limit, which is reached for about 20 modes, **2)** The decay rate $N_{tot}^{-6.5}$ applies to all cases (as expected in view of the corresponding results for $\varphi_{-2}$, Figure 5 of the main paper), **3)** The accuracy of Eq. (5.5) is only slightly affected by the reduction of smoothness of $\eta$ from $C^\infty$ to $C^4$, **4)** To achieve an accuracy $10^{-5}$ (resp., $10^{-4}$), which is satisfactory in most applications, we need six (resp., four) modes for $\varepsilon \leq 0.5$, or seven (resp., five) modes, for any $\varepsilon \leq 0.9$.

In order to assess the performance of the present approach, let us briefly discuss some other methods of computation of the DtN operator. The most popular and well-studied one is based on an operator-Taylor (perturbation) series expansion of the form

$$\mathcal{G}[\eta, h]\psi \approx \sum_{\tilde{n}=0}^{\tilde{n}=\tilde{N}} \mathcal{G}_{\tilde{n}}[\eta, h]\psi, \tag{E1}$$

where $\mathcal{G}_{\tilde{n}}[\eta, h]$ is the $\tilde{n}^{\text{th}}$-order Fréchet derivative of $\mathcal{G}[\eta, h]$ with respect to the field $\eta(\cdot)$. This approach was initially introduced in the flat-bottom case ($h = h_0 = \text{cst}$) [4] and further studied and extended by several authors; see e.g. [3], [5], [6]. The perturbation series (E1), avoiding the direct solution of Eqs. (5.1), allows for efficient computations limited, however, to cases with horizontal periodicity, flat bottom surface, and low-to-moderate slope of the upper boundary. Let it be noted that in the case where both (upper and lower) boundaries are non-planar, the expansion of $\mathcal{G}[\eta, h]$ becomes very cumbersome and difficult to implement. In contrast with Eq. (E1), the modal characterization (5.5) is not of perturbative character and, thus, it is not limited to mildly deformed domains $D_h^\eta(X)$. Further, it does not require any horizontal spatial periodicity([15]), and can easily accommodate all physically relevant lateral boundary conditions (see Appendix C). Computation of the DtN operator for the same benchmark problems, have been reported, among others, by Nicholls and Reitich [3], using three variants of the operator series method based on Eq. (E1): the operator expansion (OE) method, the field expansion (FE) method, and the transformed field expansion (TFE) method. A first qualitative comparison between these pertubative methods and the present Consistent Coupled Mode Theory (CCMT) is summarised by the following points: **1)** The OE and the FE methods always diverge as the order of truncation of the operator series exceeds some value (about 10 or smaller, dependent on $\varepsilon$ and $N_x$; see, e.g., Figs. 1 to 6 in [3]), while CCMT does not diverge for any number of terms ([16]). **2)** CCMT with 6 or 7 modes provides better accuracy in comparison with the best accuracy possible by OE and FE before their divergence. **3)** The TFE does not diverge, reaching a plateau limit, with a limiting accuracy which is slightly better than that of CCMT for small values of $\varepsilon$ (see Figs. 1 to 5 in [3] in comparison with **Figure E1**, above), but it is worse for high values of $\varepsilon$; for example, for the "rough" profile and $\varepsilon = 0.8$, the limiting accuracy of TFE is a little more than $10^{-4}$ (see Fig. 6 of [3]), while the limiting accuracy of CCMT is about $10^{-6}$. It is recalled here that the TFE is implemented by using both operator-Taylor expansion and vertical discretization of the domain in order to calculate each term of the expansion, which makes it comparable with direct numerical methods [7], [8], [9] concerning the computational burden. In contrast, the present method offers the high accuracy of direct numerical methods with substantially less computational cost, through a genuine dimensional reduction of the problem. In addition, it applies to general domains $D_h^\eta$ with both upper and lower boundaries nonplanar, and it is not limited by horizontal periodicity. An example of the calculation of the DtN operator for a strongly undulated domain (with both upper and lower boundaries non-planar) is presented below.

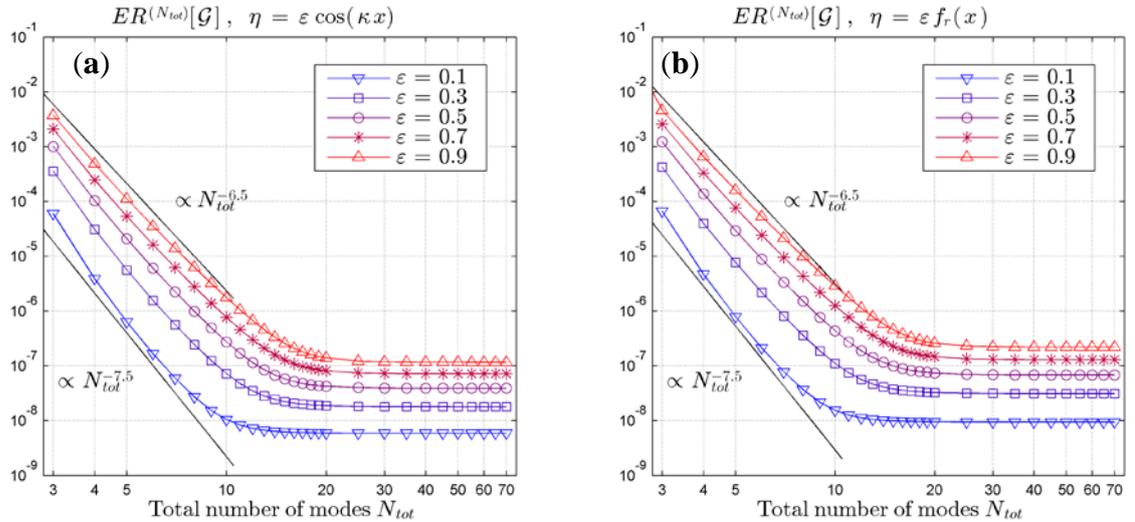

**Figure E1**. $L^2$−error of the DtN operator vs. the total number of modes used, for $\varepsilon = 0.1(0.2)0.9$. (**a**) Smooth case, (**b**) Rough case.

**(b) Computation of the DtN operator for a strip-like domain having both upper and lower boundaries non-planar**

The main advantage of the CCMT seems to us to be its ability to treat cases with (strongly) varying bathymetry $h = h(x)$ and free-surface elevation $\eta = \eta(x)$, and any type of lateral boundary conditions, without any additional cost apart from the inclusion of the bottom boundary mode, $n = -1$, in the series expansion (2.2). As a demonstration of this ability, we consider the DtN operator corresponding to the highly nonuniform domain and Dirichlet data $\psi$ shown in **Figure E2(a)**. In this domain, both upper and lower boundary surfaces are strongly "undulated", ending to a flat part with different left and right depths. Homogeneous Neumann boundary conditions are applied on the lower surface and the lateral boundaries. For comparison purposes, the DtN is computed using Eq. (5.5) with a total number of four and seven modes, and by a standard Boundary Element Method (BEM). In both methods, 1000 horizontal grid points were used. Results are shown in **Figure E2(b).** The differences between $\mathcal{G}^{(4)}[\eta, h]\psi$ and $\mathcal{G}^{(7)}[\eta, h]\psi$ are minor, and mainly restricted near the ends of the boundary nonuniformities. The BEM and CCMT results are indistinguishable in the scale of the figure. The $L^2-$ relative discrepancy between the BEM solution and $\mathcal{G}^{(4)}[\eta,h]\psi$ is of order $O(10^{-4})$. The calculations using CCMT are $10^2$ times faster, verifying the efficiency of CCMT against direct numerical methods.

---

([15]) The only reason for which we have considered herein a spatially periodic example is to compare our results with those of other authors, who compute the DtN operator using Eq. (E1).

([16]) Recall, however, that the meaning of *terms* in the two expansions is quite different. The former is an expansion with respect to vertical basis functions, the latter is an operator-Taylor expansion with respect to the horizontal field $\eta(\cdot)$.

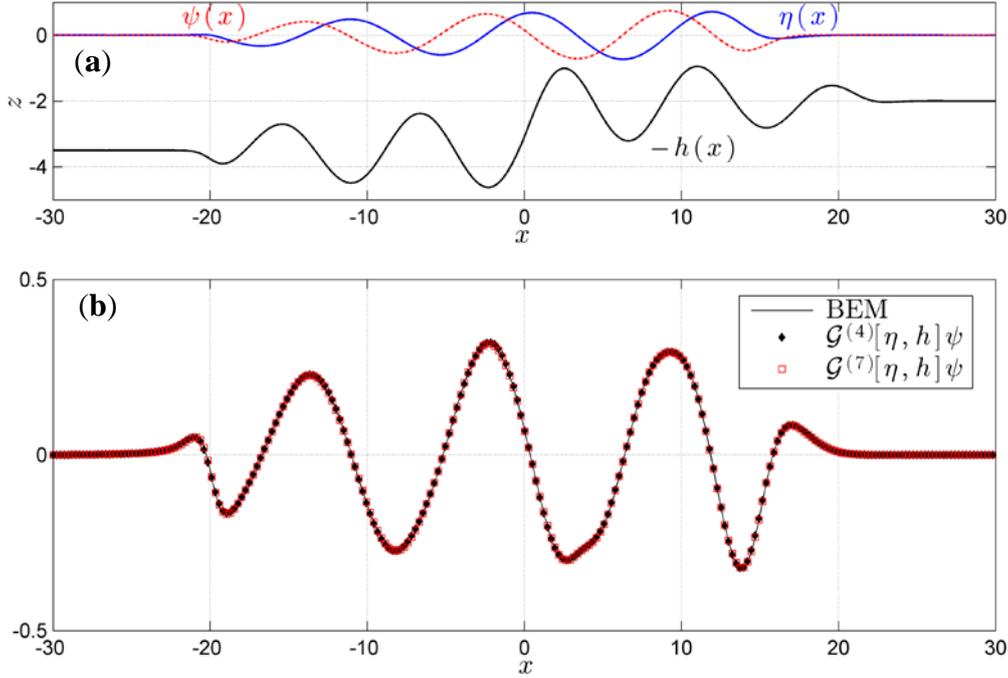

**Figure E2**: (**a**): Non-uniform domain and Dirichlet data $\psi$ used in the computations. (**b**) The DtN operator computed using CCMT, $\mathcal{G}^{(4)}[\eta,h]\psi$ and $\mathcal{G}^{(7)}[\eta,h]\psi$, and BEM.